\documentclass[11pt, reqno]{amsart}

\newif\ifrevision
\revisiontrue

\usepackage{float}

\usepackage[scaled]{helvet}
\usepackage{mathtools}
\usepackage{etoolbox}
\usepackage[utf8]{inputenc}
\usepackage[T1]{fontenc}

\usepackage{amsmath,amssymb,amsfonts,mathrsfs,amsbsy,amsthm}
\usepackage{relsize}
\usepackage{enumitem}
\usepackage{anyfontsize}
\usepackage{tikz,xcolor}
\usepackage{bm}
\usepackage[top=2.5cm,bottom=2.5cm,left=3cm,right=2.5cm]{geometry}
\usepackage[french,english]{babel}
\definecolor{lava}{rgb}{0.81,0.06,0.13}
\definecolor{Cblue}{rgb}{0.50,0.85,0.85}
\definecolor{lime}{HTML}{A6CE39}
\definecolor{linkcol}{rgb}{0.0,0.42,0.30}
\usepackage[colorlinks=true,linkcolor=lava,citecolor=blue,urlcolor=brown]{hyperref}
\definecolor{lava}{rgb}{0.8,0.2,0.0}

\definecolor{chgcol}{rgb}{0.0,0.25,0.85}  
\definecolor{chkcol}{rgb}{0.78,0.0,0.55}   
\ifrevision
  
  \newcommand{\chk}[1]{{\color{chkcol}#1}}
  \newenvironment{chgblock}{\color{chgcol}}{}
  \newenvironment{chkblock}{\color{chkcol}}{}
\else
  
  \newcommand{\chk}[1]{#1}

\fi

\DeclareRobustCommand{\orcidicon}{%
  \begin{tikzpicture}
    \draw[lime,fill=lime] (0,0) circle [radius=0.16]
      node[white] {{\fontfamily{phv}\selectfont\tiny ID}};
    \draw[white,fill=white] (-0.0625,0.095) circle [radius=0.007];
  \end{tikzpicture}\hspace{-2mm}}

\foreach \x in {A,B}{%
  \expandafter\xdef\csname orcid\x\endcsname{%
    \noexpand\href{https://orcid.org/\csname orcidauthor\x\endcsname}{\noexpand\orcidicon}}}

\numberwithin{equation}{section}

\newtheorem{theo}{Theorem}[section]
\newtheorem{prop}{Proposition}[section]
\newtheorem{lem}{Lemma}[section]
\newtheorem{coro}{Corollary}[section]
\newtheorem{conj}{Conjecture}[section]

\theoremstyle{definition}
\newtheorem{defn}{Definition}[section]
\newtheorem{example}{Example}[section]
\newtheorem{conv}{Convention}[section]

\theoremstyle{remark}
\newtheorem{rem}{Remark}[section]

\renewcommand{\geq}{\geqslant}\renewcommand{\ge}{\geqslant}
\renewcommand{\leq}{\leqslant}\renewcommand{\le}{\leqslant}

\newcommand{\R}{\mathbb{R}}

\newcommand{\RD}{\mathbb{R}^N}
\newcommand{\Lap}{(-\Delta)^{\theta/2}}

\newcommand{\CapD}[1]{\partial_t^{#1}}

\definecolor{bluegreen}{rgb}{0.0, 0.3, 0.9}

\title[Critical Exponent for a Time--Space Fractional Reaction--Diffusion Equation]
{On the Critical Exponent for a Time--Space Fractional
Reaction--Diffusion Equation with a Slowly Varying
Time-Dependent Source}

\author[Ahmad Z.~Fino \& M.~Majdoub]
{Ahmad Z.~Fino\orcidB{}\,\; \& \,\;Mohamed Majdoub\orcidA{}}
\address[Ahmad Z.~Fino]{Department of Mathematics, Faculty of Sciences, Lebanese University, P.O. Box 1352 Tripoli, Lebanon.}
\address[M.~Majdoub]{%
  Department of Mathematics, College of Science, Imam Abdulrahman Bin Faisal
  University, P.\,O.\ Box 1982, Dammam, Saudi Arabia.\newline
  Basic and Applied Scientific Research Center, Imam Abdulrahman Bin Faisal
  University, P.\,O.\ Box 1982, 31441, Dammam, Saudi Arabia.}
\email{\tt afino@ul.edu.lb}
\email{\tt mmajdoub@iau.edu.sa}



\begin{document}

\begin{abstract}
We investigate the Fujita critical exponent for the semilinear
parabolic Cauchy problem
$$
\CapD{\alpha} u +\Lap u \;=\; h(t)\, u^p, \qquad
(t,x) \in (0,\infty)\times \RD,
\qquad u(0,\cdot)=u_0\ge 0,
$$
where $\alpha\in(0,1]$, $\theta\in(0,2]$, $p>1$, $\CapD{\alpha}$ is the Caputo
derivative of order $\alpha$, and $h(t)\sim t^{\gamma}\,\ell(t)$ for large
$t$, with $\gamma\in\R$ and $\ell$ a positive continuous function that is
slowly varying at infinity in the sense of Karamata.
For the classical time derivative $\alpha=1$ and every $\theta\in(0,2]$
we obtain a complete classification. If $\gamma\le-1$, small data give
global solutions for every $p>1$, so there is no Fujita threshold. If
$\gamma>-1$, the critical exponent is $p_c(\gamma)=1+\theta(\gamma+1)/N$:
every nontrivial nonnegative solution blows up when $1<p<p_c(\gamma)$,
while small global solutions exist when $p>p_c(\gamma)$. At
$p=p_c(\gamma)$ the slowly varying factor decides. Blow-up occurs when
$\int_1^{\infty} t^{-1}\ell(t)\,dt = \infty$, whereas small global
solutions exist when this integral converges. Hence global existence may
occur \emph{at} the critical exponent, in contrast with the autonomous case.
{For a genuinely fractional time derivative $\alpha\in(0,1)$, the
memory of the Caputo derivative replaces $\gamma+1$ by $\gamma+\alpha$.
We prove small-data global existence for every $p>1$ when
$\gamma\le-\alpha$. When $\gamma>-\alpha$, we prove in every dimension
that all nontrivial solutions blow up for
$1<p<p_c(\gamma,\alpha):=1+\theta(\gamma+\alpha)/(N\alpha)$. If $N<\theta$
(that is, $N=1$ and $\theta\in(1,2]$), we also obtain global existence for
$p>p_c(\gamma,\alpha)$, so the classification is complete; at
$p=p_c(\gamma,\alpha)$ the dichotomy is then governed by the boundedness of
$\ell$ rather than by the Karamata integral. If $N\ge\theta$, we obtain
global existence for $p>1+(\gamma+\alpha)/\alpha$, a threshold that
coincides with $p_c(\gamma,\alpha)$ when $N=\theta$. Finally, we exhibit a
second, purely local blow-up mechanism which shows that, for $N>\theta$
and $\gamma>\alpha\theta/(N-\theta)$, the exponent $p_c(\gamma,\alpha)$ is
\emph{not} the Fujita exponent.}
Throughout, the source is understood as a large-time profile completed
near $t=0$ so as to be continuous up to the origin
(Convention~\ref{conv:source}). This normalisation does not affect any of
the critical exponents above, but it is needed for the Cauchy problem to
be well posed.
Our results extend the classical theorems of Meier [Arch.~Rational
Mech.~Anal. {\bf 109} (1990)] to the nonlocal setting and clarify the role
of Karamata's slow variation at the critical exponent.
\end{abstract}
\maketitle

\section{Introduction}
\label{sec:intro}

The long-time behavior of solutions to semilinear parabolic equations of
the form
\begin{equation}
\label{eq:main}
\CapD{\alpha} u +\Lap u \;=\; h(t)\, u^p,
\qquad h(t)\sim t^{\gamma}\,\ell(t),
\qquad (t,x)\in(0,\infty)\times\RD,
\end{equation}
with nonnegative initial data $u(0,\cdot)=u_0\ge 0$, $u_0\not\equiv 0$,
has attracted considerable interest since the pioneering work of
Fujita~\cite{Fujita66}, which was motivated in part by the earlier
eigenfunction-based blow-up criterion of Kaplan~\cite{Kaplan63}.
For the autonomous classical case ($h\equiv 1$, $\alpha=1$, $\theta=2$),
Fujita identified the remarkable threshold
$$
p_F \;=\; 1 + \tfrac{2}{N}.
$$
If $1<p<p_F$, every nonnegative nontrivial solution blows up in finite
time. If $p>p_F$, global positive solutions exist for sufficiently small
initial data. The borderline case $p=p_F$ belongs to the blow-up regime.
This was shown by Hayakawa~\cite{Hayakawa73} for $N=1,2$ and by
Kobayashi, Sirao and Tanaka~\cite{KST77} for general $N$, and
Weissler~\cite{Weissler81} later gave a much shorter proof valid in any
dimension; see also Aronson and Weinberger~\cite{AronsonWeinberger78}.
This phenomenon, now called the \emph{Fujita phenomenon}, is robust under
many perturbations. It has been extended to systems, to nonlocal
diffusion and to various geometries; we refer to the surveys
{\cite{Levine90,DengLevine00} and to the monograph
\cite{QuittnerSouplet}} for an overview.

A natural and physically relevant generalisation concerns time-dependent
source intensities. When $h(t)$ behaves like a power $t^q$, or more
generally like a regularly varying function, the critical exponent shifts
in a predictable way. For the heat equation on a domain $D\subset\RD$
with Dirichlet boundary conditions, this question was treated in great
generality by Meier~\cite{Meier90}, the case $D=\RD$ being included; see
also his earlier works~\cite{Meier86,Meier88}.
{Coefficients depending on both $t$ and $x$ were considered, for
instance, by Qi~\cite{Qi98} and Igarashi~\cite{Igarashi2004}, while
Pinsky~\cite{Pinsky97} analysed the spatially inhomogeneous problem
$u_t=\Delta u+a(x)u^p$ in $\RD$. For the nonlocal equation
$u_t+\Lap u=h(t)u^p$, nonexistence of global solutions was established by
Guedda and Kirane~\cite{GueddaKirane}, extending the work of
Sugitani~\cite{Sugitani75}.}
Meier's fundamental observation is that the critical exponent admits the
universal formula
\begin{equation}\label{eq:Meier-formula}
\begin{split}
p^{\,\ast} \;&=\; 1 + \frac{q+1}{s^{\,\ast}},\\
s^{\,\ast}\;&:=\;\sup\bigl\{s>0: \exists\, w\not\equiv 0:\; \partial_t w=\Delta w\,\; \text{with } \limsup_{t\to\infty} t^s\|w(t)\|_{\infty}<\infty\bigr\},
\end{split}
\end{equation}
where the supremum runs over nonnegative solutions of the homogeneous
problem. For $D=\RD$ one has $s^{\,\ast}=N/2$, so that
$p^{\,\ast}=1+2(q+1)/N${, which reduces to Fujita's value $p_F$ when
$q=0$}. On conical or sectorial domains, $s^{\,\ast}$ encodes the
geometry through the first Dirichlet eigenvalue of the
Laplace--Beltrami operator on the spherical cap;
see~\cite{Meier90,BandleLevine89,LevineMeier90}.

{A closely related line of research concerns \emph{additive} forcing
terms, as in $u_t-\Delta u=|u|^p+t^{\sigma}w(x)$, for which the Fujita
exponent is again shifted by the forcing; see Bandle, Levine and
Zhang~\cite{BandleLevineZhang00}, the second author's
works~\cite{Majdoub23,AlshehriEtAl23} on Hardy--H\'enon and forced
parabolic equations, and~\cite{BenBelgacemMajdoub26}, where a Caputo
time derivative and a fractional Laplacian are combined with a forcing
term $t^{\sigma}w(x)$. Sharp lifespan estimates for fractional
Hardy--H\'enon parabolic equations are obtained
in~\cite{MajdoubTorebek26}. The present work complements these papers:
the time dependence now enters multiplicatively through the
intensity~$h(t)$ of the nonlinear source.}

The purpose of the present paper is twofold.
\begin{itemize}
\item[(a)] We \emph{lift} the Meier framework to the nonlocal setting.
The Laplacian is replaced by the fractional Laplacian $\Lap$ of order
$\theta\in(0,2)$ and, more delicately, the time derivative is replaced
by the Caputo derivative $\CapD{\alpha}$ of order $\alpha\in(0,1)$. For
the space-fractional equation with classical time derivative
($\alpha=1$, $\theta\in(0,2)$), the autonomous Fujita threshold is
$p_F=1+\theta/N$, and blow-up holds up to and including the borderline;
this goes back to Sugitani~\cite{Sugitani75}. For the time-fractional
equation ($\alpha\in(0,1)$), the autonomous threshold is \emph{again}
$p_F=1+\theta/N$, independently of $\alpha$. Now, however, the borderline
$p=p_F$ admits global solutions; see Zhang and Sun~\cite{Zhang14} for
$\theta=2$, and Cort\'azar, Quir\'os and
Wolanski~\cite{CortazarQuirosWolanski24} for the fully nonlocal case.
The memory of the Caputo derivative thus affects the role of $\alpha$ in
a subtle way. It does not shift the autonomous exponent but, as we shall
see, it \emph{does} shift the exponent of a time-dependent source,
replacing Meier's numerator $\gamma+1$ by $\gamma+\alpha$. Heuristically,
this is the trace of the kernel $(t-s)^{\alpha-1}$ in the Duhamel
representation; see \S\ref{subsec:meier-frac}.
{Moreover, when $N>\theta$, a second and purely local mechanism
appears. It is caused by the slow decay $t^{-\alpha}$ of the fundamental
solution near the origin and produces blow-up for $p<1+\gamma/\alpha$,
a range that may exceed $p_c(\gamma,\alpha)$; see
Proposition~\ref{prop:local}.}

\item[(b)] We \emph{incorporate} a slowly varying prefactor $\ell$ in
the source, that is, a measurable function
$\ell:(0,\infty)\to(0,\infty)$ such that $\ell(\lambda t)/\ell(t)\to 1$
as $t\to\infty$ for every $\lambda>0$. Such functions appear naturally in
physics whenever scaling symmetries are only approximate, or when
logarithmic corrections arise from boundary layers, mean-field limits or
anomalous transport in disordered media \cite{Bingham87,Resnick07}. They
strictly include the positive constants, and also logarithms, iterated
logarithms and powers of logarithms. Our analysis shows that the slowly
varying factor never shifts the \emph{value} of the critical exponent.
At $\alpha=1$, however, it plays a decisive role \emph{at} the critical
exponent $p=p_c(\gamma)$, for every $\gamma>-1$: the dichotomy between
blow-up and small-data global existence is determined by the convergence
or divergence of $\int^{\infty} t^{-1}\ell(t)\,dt$. {For
$\alpha\in(0,1)$ and $N<\theta$, the same dichotomy is instead governed by
the boundedness of $\ell$.}
Throughout, the source $h$ is normalised near the origin so as to be
continuous up to $t=0$ (Convention~\ref{conv:source}). This is
immaterial for every critical exponent below, but it is needed for the
Cauchy problem to be well posed. Indeed, the literal density
$t\mapsto t^{\gamma}\ell(t)$ need not be integrable at the origin (with
$\ell\equiv 1$, it is not integrable as soon as $\gamma\le-1$), and slow
variation constrains only the regime $t\to\infty$.
\end{itemize}

\noindent
\textbf{Main results (informal).} For $\alpha=1$ and $\theta\in(0,2]$ we
obtain the complete classification
$$
p_c(\gamma)\;=\;
\begin{cases}
\infty,
& \text{if }\gamma\le-1,\\[6pt]
1+\dfrac{\theta(\gamma+1)}{N},
& \text{if }\gamma>-1.
\end{cases}
$$
Moreover, for every $\gamma>-1$, the behavior \emph{at}
$p=p_c(\gamma)$ is governed by the integrability of
$t\mapsto t^{-1}\ell(t)$ at infinity. Divergence forces blow-up of every
nontrivial solution, whereas convergence allows small-data global
existence.
{For $\alpha\in(0,1)$ we prove small-data global existence for every
$p>1$ when $\gamma\le-\alpha$. When $\gamma>-\alpha$, we prove blow-up of
every nontrivial solution for $1<p<p_c(\gamma,\alpha):=
1+\theta(\gamma+\alpha)/(N\alpha)$ in every dimension, and also at
$p=p_c(\gamma,\alpha)$ when $\ell$ is unbounded. When $N<\theta$, i.e.\
$N=1$ and $\theta\in(1,2]$, small global solutions exist for
$p>p_c(\gamma,\alpha)$, and also at $p=p_c(\gamma,\alpha)$ when $\ell$ is
bounded; the classification is therefore complete. When $N\ge\theta$, we
reach the threshold $1+(\gamma+\alpha)/\alpha$, and we show that, for
$N>\theta$, blow-up also occurs for $1<p<1+\gamma/\alpha$; we conjecture
that the Fujita exponent is then $\max\{p_c(\gamma,\alpha),
1+\gamma/\alpha\}$ (Conjecture~\ref{conj:frac}).}
Note that $p_c(\gamma,1)=p_c(\gamma)$ and $p_c(0,\alpha)=1+\theta/N$ for
every $\alpha$, in agreement with~\cite{Zhang14,CortazarQuirosWolanski24}.

\medskip
The paper is organised as follows.
Section~\ref{sec:prel} collects the necessary definitions and classical
results: slowly varying functions (\S\ref{subsec:slow}), the Caputo
derivative and the fractional Laplacian
(\S\S\ref{subsec:caputo}--\ref{subsec:flap}), the time--space fractional
solution operators and their decay (\S\ref{subsec:semigroup}), mild
solutions and the standing convention on the source
(\S\ref{subsec:mild}), the classical Fujita theory
(\S\ref{subsec:fujita}) and Meier's theorem (\S\ref{subsec:meier}).
Section~\ref{sec:main} contains the precise statements of the main
results for $\alpha=1$ and for $\alpha<1$, the accompanying conjecture,
and a discussion of their consequences.
Section~\ref{sec:other} surveys related classes of $h(t)$ (exponential,
integrable, oscillatory) and their Fujita exponents.
{Section~\ref{sec:methods} contains the proofs: a fixed-point argument
in weighted spaces for global existence, a Kaplan-type functional for
blow-up at $\alpha=1$, a nonlinear iteration of
Cort\'azar--Quir\'os--Wolanski type for blow-up when $\alpha<1$, and a
local iteration exhibiting the second blow-up mechanism.}
We close with a few open problems in Section~\ref{sec:conc}.

\section{Preliminaries}
\label{sec:prel}

\subsection{Slowly varying functions}
\label{subsec:slow}

\begin{defn}
A measurable function $\ell:(0,\infty)\to(0,\infty)$ is called
\emph{slowly varying at infinity} if
\begin{equation}\label{eq:sv-def}
\lim_{t\to\infty}\frac{\ell(\lambda t)}{\ell(t)}=1
\qquad\text{for every }\ \lambda>0.
\end{equation}
A measurable function $g:(0,\infty)\to(0,\infty)$ is \emph{regularly
varying of index $\rho\in\R$} if $g(t)=t^\rho \ell(t)$ for some slowly
varying $\ell$.
\end{defn}

Typical examples are
$$
\ell(t) = \bigl(\log(1+t)\bigr)^a,\quad
\ell(t) = \bigl(\log\log(2+t)\bigr)^b,\quad
\ell(t) = \exp\!\bigl((\log t)^\delta\bigr)\;\; (0<\delta<1,\ t\ge 1),
$$
as well as products and quotients of such functions. Any function that
converges to a positive constant at infinity is trivially slowly varying.

By Karamata's representation theorem \cite[Theorem 1.3.1]{Bingham87}, a
slowly varying function can be written as
\begin{equation}\label{eq:karamata-rep}
\ell(t)\;=\;c(t)\,\exp\!\left(\int_{a}^{t}\frac{\varepsilon(s)}{s}\,ds\right),
\qquad t\ge a,
\end{equation}
for some $a>0$, where $c(t)\to c_0\in(0,\infty)$ and $\varepsilon(t)\to 0$
as $t\to\infty$. We shall repeatedly use the following four
consequences.

\begin{lem}\label{lem:karamata}
Let $\ell$ be slowly varying at infinity {and locally bounded on
$[a,\infty)$ for some $a>0$ (this holds if $\ell$ is continuous)}. Then
the following properties hold.
\begin{enumerate}
\item[(K1)] \emph{(Sub-polynomial growth and decay.)} For every
$\varepsilon>0$, $t^{-\varepsilon}\ell(t)\to 0$ and
$t^{\varepsilon}\ell(t)\to\infty$ as $t\to\infty$
{\cite[Proposition~1.3.6]{Bingham87}}.
\item[(K2)] \emph{(Karamata's integral theorem.)} For $\sigma>-1$,
$\displaystyle\int_a^{t} s^{\sigma}\ell(s)\,ds\sim \frac{t^{\sigma+1}\ell(t)}{\sigma+1}$
as $t\to\infty$. For $\sigma<-1$, we have
$\int_a^{\infty} s^{\sigma}\ell(s)\,ds<\infty$ and
$\displaystyle \int_t^{\infty} s^{\sigma}\ell(s)\,ds\sim \frac{t^{\sigma+1}\ell(t)}{-(\sigma+1)}$
{\cite[Propositions~1.5.8 and~1.5.10]{Bingham87}}.
Both statements involve a fixed lower limit $a>0$. They describe the
tail $t\to\infty$ and say nothing about integrability at the origin.
\item[(K3)] \emph{(Borderline case $\sigma=-1$.)} The integral
$\int^{\infty} s^{-1}\ell(s)\,ds$ may converge or diverge, depending on
$\ell$. For example, $\ell\equiv 1$ gives divergence, $\ell(t)=(\log t)^{-2}$
gives convergence, and $\ell(t)=(\log t)^{-1}$ gives divergence (like
$\log\log t$).
\item[(K4)] \emph{(Decay under a convergent borderline integral.)} If
$\int^{\infty} s^{-1}\ell(s)\,ds<\infty$, then $\ell(t)\to 0$ as
$t\to\infty$. Indeed, by the Uniform Convergence Theorem
\cite[Theorem 1.2.1]{Bingham87}, $\ell(s)\ge (1+o(1))\,\ell(t)$
uniformly for $s\in[t,2t]$. We deduce that
$\int_t^{2t}s^{-1}\ell(s)\,ds \ge (\log 2)(1+o(1))\,\ell(t)$, and the
left-hand side tends to $0$ as $t\to\infty$ since the tail integral
converges.
\end{enumerate}
\end{lem}

{We shall also use the Uniform Convergence Theorem in the following
form: for every compact interval $[\lambda_1,\lambda_2]\subset(0,\infty)$,
\begin{equation}\label{eq:UCT}
\sup_{\lambda\in[\lambda_1,\lambda_2]}\Bigl|\frac{\ell(\lambda t)}{\ell(t)}-1\Bigr|
\longrightarrow 0\qquad (t\to\infty),
\end{equation}
see \cite[Theorem 1.2.1]{Bingham87}.}

\subsection{Caputo fractional derivative}
\label{subsec:caputo}

\begin{defn}
Let $\alpha\in(0,1)$. For an absolutely continuous function
$f:[0,T]\to\R$, the Caputo derivative of order $\alpha$ is defined by
\begin{equation}\label{eq:caputo}
\CapD{\alpha} f(t)
\;=\;\frac{1}{\Gamma(1-\alpha)}\int_0^{t}(t-s)^{-\alpha}\,f'(s)\,ds,
\qquad t>0.
\end{equation}
For $\alpha=1$ we set $\CapD{1}f=f'$ in the classical sense.
\end{defn}

The Caputo derivative is the natural notion of fractional time
derivative for Cauchy problems with classical initial conditions. Since
$\CapD{\alpha}c=0$ for every constant $c$, the initial value is
consistent with the equation. Equivalently, on absolutely continuous
functions, $\CapD{\alpha}f(t) = D^{\alpha}\!\bigl[f(t)-f(0)\bigr]$, where
$D^{\alpha}$ is the Riemann--Liouville derivative;
see~\cite{KilbasBook,Podlubny99}.

\subsection{Fractional Laplacian}
\label{subsec:flap}

For $\theta\in(0,2)$, the fractional Laplacian $\Lap$ acts on Schwartz
functions $\varphi\in\mathcal{S}(\RD)$ by
\begin{equation}\label{eq:flap-fourier}
\widehat{\Lap\varphi}(\xi)\;=\;|\xi|^{\theta}\,\widehat{\varphi}(\xi),
\end{equation}
or, equivalently, by the singular integral
\begin{equation}\label{eq:flap-singular}
\Lap\varphi(x)
\;=\;C_{N,\theta}\,
\mathrm{P.V.}\!\int_{\RD}\frac{\varphi(x)-\varphi(y)}{|x-y|^{N+\theta}}\,dy,
\end{equation}
with the explicit constant
$C_{N,\theta}=\theta\,2^{\theta-1}\Gamma((N+\theta)/2)/\bigl(\pi^{N/2}\Gamma(1-\theta/2)\bigr)$;
see~{\cite{DiNezzaPV12,DipierroValdinoci17,KwasnickiTen17}}.
The operator extends to a nonnegative self-adjoint operator on
$L^{2}(\RD)$ with form domain $H^{\theta/2}(\RD)$. For $\theta=2$ we
recover $-\Delta$.

\subsection{The time--space fractional solution operators and kernel estimates}
\label{subsec:semigroup}

Consider the homogeneous Cauchy problem
\begin{equation}\label{eq:hom}
\CapD{\alpha} w + \Lap w = 0,\qquad w(0,\cdot)=w_0.
\end{equation}
Using the Fourier transform in $x$ and the Laplace transform in $t$, we
obtain $w(t,x)= [S_{\alpha,\theta}(t)w_0](x)$, where
\begin{equation}\label{eq:S-operator}
S_{\alpha,\theta}(t)w_0(x)
\;=\;\int_{\RD} Z(t,x-y)\,w_0(y)\,dy,
\end{equation}
and the fundamental solution $Z$ is characterised by
$\widehat{Z}(t,\xi)=E_\alpha(-t^{\alpha}|\xi|^{\theta})$. Here $E_\alpha$
is the Mittag-Leffler function
$E_\alpha(z)=\displaystyle\sum_{k=0}^{\infty}\frac{z^k}{\Gamma(\alpha k+1)}$;
see \cite{MainardiLuchkoPagnini01,KimLim15,KSZ17}. For every $t>0$ the
kernel $Z(t,\cdot)$ is a positive probability density, so that
$\int_{\RD}Z(t,x)\,dx=1$. {Since $\widehat{Z}(t,\xi)$ depends only on
$t^{\alpha}|\xi|^{\theta}$, the kernel is self-similar:
\begin{equation}\label{selfsimilarZ}
Z(t,x)=t^{-\frac{N\alpha}{\theta}}F\bigl(xt^{-\frac{\alpha}{\theta}}\bigr),
\qquad (t,x)\in(0,\infty)\times\RD,
\end{equation}
where the profile $F=Z(1,\cdot)$ is radial, positive and continuous on
$\RD\setminus\{0\}$. Writing $f\asymp g$ when $c^{-1}g\le f\le cg$ for
some constant $c\ge1$, the two-sided estimates of
\cite[Lemma~3.3]{KSZ17}  give, for
$0<|y|\le1$,
$$F(y)\asymp\left\{\begin{array}{ll}
|y|^{-N+\theta}&\hbox{if $N>\theta$ and $0<\alpha<1$},\\[2pt]
1+|\log(|y|^{\theta})|&\hbox{if $N=\theta$ and $0<\alpha<1$},\\[2pt]
1&\hbox{if $\alpha=1$, or if $N<\theta$ and $0<\alpha<1$}.
\end{array}
\right.
$$}
In particular, {since $F$ is bounded from below by a positive constant
on the unit ball in all three cases,}
\begin{equation}\label{selfsimilarZZ}
Z(t,x)\geq C\,t^{-\frac{N\alpha}{\theta}},\qquad\hbox{if}\,\,|x|\leq t^{\alpha/\theta},
\end{equation}
{and, if $N>\theta$ and $0<\alpha<1$,
\begin{equation}\label{eq:Z-sing}
Z(t,x)\geq C\,t^{-\alpha}|x|^{\theta-N},\qquad\hbox{if}\,\,0<|x|\leq t^{\alpha/\theta}.
\end{equation}}
For $\alpha=1$ and $\theta=2$ we recover the Gaussian kernel, and for
$\alpha=1$, $\theta\in(0,2)$, the $\theta$-stable kernel.

The validity of $L^p$--$L^q$ smoothing estimates depends crucially on
$\alpha$. The following statement records the sharp ranges; see
\cite{KimLim15,KSZ17,CortazarQuirosWolanski24}.

\begin{prop}\label{prop:heat-kernel}
Let $\theta\in(0,2]$ and $1\le p\le q\le\infty$.
\begin{enumerate}
\item[\textup{(i)}] If $\alpha=1$, then for all $\phi\in L^p(\RD)$,
\begin{equation}\label{eq:Lp-Lq}
\|S_{\alpha,\theta}(t)\phi\|_{L^q(\RD)}
\;\le\; C\, t^{-\frac{N\alpha}{\theta}\bigl(\frac1p-\frac1q\bigr)}\,
\|\phi\|_{L^p(\RD)},\qquad t>0,
\end{equation}
with no restriction on $(p,q)$.
\item[\textup{(ii)}] If $\alpha\in(0,1)$, then~\eqref{eq:Lp-Lq} holds
whenever $\tfrac1p-\tfrac1q<\tfrac{\theta}{N}$ and fails whenever
$\tfrac1p-\tfrac1q>\tfrac{\theta}{N}$.
\end{enumerate}
When it is valid, the decay rate in~\eqref{eq:Lp-Lq} is sharp. In both
cases $S_{\alpha,\theta}(t)$ is a contraction on every $L^r(\RD)$,
$1\le r\le\infty$.
\end{prop}

{Part (ii) follows from Young's inequality, since by
\eqref{selfsimilarZ} we have
$\|Z(t,\cdot)\|_{L^m}=t^{-\frac{N\alpha}{\theta}(1-\frac1m)}\|F\|_{L^m}$,
and $F\in L^m(\RD)$ if and only if $(N-\theta)m<N$, that is,
$1-\frac1m<\frac{\theta}{N}$.}

\begin{rem}\rm\label{rem:kernel-singularity}
\leavevmode
\begin{enumerate}
\item[\textup{(i)}] For $\alpha\in(0,1)$, consider the critical line
$\tfrac1p-\tfrac1q=\tfrac{\theta}{N}$, which is non-void only when
$N\ge\theta$. There, the estimate~\eqref{eq:Lp-Lq} holds for
$1<p<q<\infty$. This is a consequence of the global pointwise bound
$Z(t,x)\le C\,t^{-\alpha}|x|^{\theta-N}$ for $N>\theta$ and of the
Hardy--Littlewood--Sobolev inequality. It fails exactly at the endpoint
pairs $(p,q)=\bigl(1,\tfrac{N}{N-\theta}\bigr)$ and
$(p,q)=\bigl(\tfrac{N}{\theta},\infty\bigr)$ (at $(p,q)=(1,\infty)$ when
$N=\theta$). In particular, \eqref{eq:Lp-Lq} holds for all
$1\le p\le q\le\infty$ when $N<\theta$, whereas the endpoint
$L^1\to L^\infty$ estimate \emph{fails} whenever $N\ge\theta$.
\item[\textup{(ii)}] The restriction in
Proposition~\ref{prop:heat-kernel}\textup{(ii)} reflects a genuine
spatial singularity of the fundamental solution. For $\alpha<1$, the
profile in \eqref{selfsimilarZ} satisfies $F(\xi)\asymp|\xi|^{\theta-N}$
as $\xi\to 0$ when $N>\theta$, with a logarithmic singularity when
$N=\theta$, so that $Z(t,\cdot)\notin L^\infty(\RD)$; see
\cite{KimLim15,KSZ17}. Consequently, for $\alpha<1$ and $N>\theta$, the
sup-norm of $S_{\alpha,\theta}(t)w_0$ with $0\le w_0\in L^1\cap L^\infty$,
$w_0\not\equiv0$, decays exactly at the rate $t^{-\alpha}$, which is
slower than the formal $L^1\to L^\infty$ rate $t^{-N\alpha/\theta}$.
{This slow decay is the origin of the local blow-up mechanism of
Proposition~\ref{prop:local}.}
The operator family $P_{\alpha,\theta}$ of \S\ref{subsec:mild} obeys
analogous estimates on its own admissible ranges{, namely with
$\theta$ replaced by $2\theta$ in the admissibility condition, because its
kernel behaves like $|x|^{2\theta-N}$ near the origin when $N>2\theta$}; see
\cite{CortazarQuirosWolanski24,KSZ17}.
\item[\textup{(iii)}] At $\alpha=1$ the exponent $N/\theta$ plays the
role of $s^{\,\ast}$ in Meier's formula~\eqref{eq:Meier-formula}. For
$\alpha<1$ the correct analogue is subtler; see \S\ref{subsec:meier-frac}.
\end{enumerate}
\end{rem}

\subsection{Mild solutions and the standing convention on the source}
\label{subsec:mild}

Following the standard approach \cite{KiraneLaskriTatar05,Zacher09}{, see also
\cite{Zhang14,CortazarQuirosWolanski24}}, the Cauchy problem~\eqref{eq:main}
is written in the integral form
\begin{equation}\label{eq:mild}
u(t,x)\;=\;S_{\alpha,\theta}(t)u_0(x)
\;+\;\int_0^{t}(t-s)^{\alpha-1}P_{\alpha,\theta}(t-s)
\bigl[h(s)\,u(s,\cdot)^p\bigr](x)\,ds,
\end{equation}
where, for each $t>0$, $P_{\alpha,\theta}(t)$ is the operator with
Fourier symbol
$\widehat{P_{\alpha,\theta}(t)f}(\xi)=\alpha\, E'_\alpha\!\bigl(-t^{\alpha}|\xi|^{\theta}\bigr)\,\widehat{f}(\xi)$.
{Since $\alpha E'_\alpha=E_{\alpha,\alpha}$, where
$E_{\alpha,\alpha}(z)=\sum_{k\ge0}z^k/\Gamma(\alpha k+\alpha)$, this
symbol equals $E_{\alpha,\alpha}(-t^{\alpha}|\xi|^{\theta})$.}
For $\alpha=1$ we have $P_{1,\theta}=S_{1,\theta}$, and~\eqref{eq:mild}
reduces to the standard Duhamel formula. The kernel of
$P_{\alpha,\theta}(t)$ is nonnegative with mass $1/\Gamma(\alpha)$.
Hence $P_{\alpha,\theta}(t)$ is bounded on every $L^r(\RD)$, uniformly in
$t$, and smoothing estimates of the form~\eqref{eq:Lp-Lq} hold on the
corresponding admissible ranges; see~\cite{CortazarQuirosWolanski24}.

\begin{defn}\label{def:mild}
Let $0\le u_0\in L^1(\RD)\cap L^\infty(\RD)$ and $T\in(0,\infty]$. A
nonnegative measurable function $u$ with
$u\in L^\infty\bigl((0,T');L^1(\RD)\cap L^\infty(\RD)\bigr)$ for every
$T'<T$ is a \emph{mild solution} of~\eqref{eq:main} on $[0,T)$ if
\eqref{eq:mild} holds for almost every $(t,x)\in(0,T)\times\RD$. For
$q\in[1,\infty]$ we say, following \cite{CortazarQuirosWolanski24}, that
$u$ is an \emph{$L^q$-solution} if
$u\in L^\infty_{\mathrm{loc}}([0,T);L^q(\RD))$ and \eqref{eq:mild} holds
in $L^q(\RD)$ for almost every $t\in(0,T)$. A mild solution in the above
sense is both an $L^1$- and an $L^\infty$-solution.
\end{defn}

\begin{conv}\label{conv:source}
The function $\ell:(0,\infty)\to(0,\infty)$ is continuous and slowly
varying at infinity, and the source intensity $h$ is understood as a
\emph{large-time profile}. Following Meier, we assume that
$$
h\in C\bigl([0,\infty)\bigr),\qquad h(t)>0\ \text{ for } t>0,
\qquad\text{and}\qquad h(t)=t^{\gamma}\ell(t)\ \text{ for } t\ge 1.
$$
The threshold $1$ is arbitrary and can be replaced by any $t_0>0$. In
particular, $h\in L^{1}_{\mathrm{loc}}\bigl([0,\infty)\bigr)$. This
normalisation is necessary for the Duhamel formulation~\eqref{eq:mild} to
make sense. Indeed, for $u$ bounded near $t=0$, the contribution of
small times $s$ to the nonlinear term behaves like
$t^{\alpha-1}\|u_0\|_\infty^{p}\int_{0}h(s)\,ds$. The literal density
$t^{\gamma}\ell(t)$, however, need not be integrable at the origin (for
$\ell\equiv 1$ it is not integrable as soon as $\gamma\le-1$), because
slow variation constrains only the behavior as $t\to\infty$. The
completion on $[0,1]$ is immaterial: every Fujita exponent below is
governed solely by the behavior of $h$ as $t\to\infty$.
For definiteness one may take $h(t)=(1+t)^{\gamma}\,\ell(1+t)$. This
function is continuous and positive on $[0,\infty)$ and
asymptotic to $t^{\gamma}\ell(t)$. {Strictly speaking, it does not
coincide with $t^\gamma\ell(t)$ for $t\ge1$, but it can be written as
$t^{\gamma}\tilde\ell(t)$ with
$\tilde\ell(t):=(1+1/t)^{\gamma}\ell(1+t)$, which is continuous,
positive and slowly varying by~\eqref{eq:UCT}. Moreover,
$\tilde\ell/\ell\to1$, so $\tilde\ell$ and $\ell$ share all the
properties used below (boundedness, unboundedness, and convergence of
$\int^\infty t^{-1}\ell(t)\,dt$).} In the subintegrable regime
$\gamma<-1$, this realisation satisfies $\int_0^\infty h(t)\,dt<\infty$.
\end{conv}

{Under Convention~\ref{conv:source}, the source is continuous on
$[0,\infty)$. A standard contraction argument in
$L^\infty((0,T);L^1\cap L^\infty)$ then provides, for every
$0\le u_0\in L^1\cap L^\infty$, a unique maximal mild solution on
$[0,T^{\ast})$, which is nonnegative. If $T^{\ast}<\infty$, then
$\limsup_{t\uparrow T^{\ast}}\|u(t)\|_{L^\infty}=\infty$; see
\cite{Weissler80} for $\alpha=1$ and \cite{Zhang14,CortazarQuirosWolanski24}
for $\alpha\in(0,1)$, whose arguments apply verbatim to a bounded
continuous coefficient $h$ on $[0,T]$.}

\subsection{Classical Fujita theory}
\label{subsec:fujita}

Consider the autonomous Cauchy problem
\begin{equation}\label{eq:autonomous}
\partial_t u - \Delta u = u^p,\qquad u(0,x)=u_0(x)\ge 0,\quad x\in\R^N.
\end{equation}
Fujita~\cite{Fujita66} proved the following dichotomy.
\begin{itemize}
\item If $1<p<1+2/N$, every nonnegative nontrivial solution blows up in
finite time.
\item If $p>1+2/N$, there exist small initial data for which the
solution exists globally and decays to zero.
\end{itemize}
The borderline exponent $p_F=1+2/N$ belongs to the blow-up case
\cite{Hayakawa73,KST77,Weissler81}.

For the autonomous nonlocal equation
$\CapD{\alpha} u + \Lap u = u^p$ in $\RD$, the threshold is
\begin{equation}\label{eq:autonomous-Fujita}
p_F\;=\;1+\frac{\theta}{N},
\end{equation}
\emph{independently of} $\alpha\in(0,1]$. The status of the borderline,
however, depends on $\alpha$.
\begin{itemize}
\item For $\alpha=1$ and $\theta\in(0,2)$, Sugitani~\cite{Sugitani75}
proved blow-up of every nontrivial nonnegative solution for
$1<p\le 1+\theta/N$, while small-data global solutions exist above the
threshold.
\item For $\alpha\in(0,1)$, Zhang and Sun~\cite{Zhang14} ($\theta=2$) and
Cort\'azar, Quir\'os and Wolanski~\cite{CortazarQuirosWolanski24}
($\theta\in(0,2)${; the details are written there for $N>2\theta$})
proved blow-up for $1<p<1+\theta/N$. In sharp contrast with the case
$\alpha=1$, they also proved the \emph{existence of global solutions at
the borderline} $p=1+\theta/N$ for suitable small data.
\end{itemize}
Test-function (nonexistence) results for related problems with
spatio-temporal fractional derivatives were obtained by Kirane, Laskri
and Tatar~\cite{KiraneLaskriTatar05}, using the method of Mitidieri and
Pohozaev~\cite{MitidieriPohozaev01}. As is typical of pure scaling
arguments, the ranges they obtain are not sharp for $\alpha<1$
(compare with \S\ref{subsec:proof-blowup-frac}).

\subsection{\texorpdfstring{Meier's theorem and its $s^{\,\ast}$-characterisation}{Meier's theorem and its s*-characterisation}}
\label{subsec:meier}

We now state Meier's results in a form suited to our adaptation.
Throughout this subsection, $\alpha=1$, $\theta=2$, and
$h\in C([0,\infty))$ with $h\ge 0$.

\begin{theo}[Meier {\cite[Theorem 2]{Meier90}}]
\label{thm:Meier2}
Consider the Cauchy problem $\partial_t u-\Delta u=h(t)u^p$ in $\RD$ with
$u(0,\cdot)=u_0\ge 0$.

\textup{(i)} \emph{(Global existence.)} If there exists a nonnegative
solution $w\not\equiv 0$ of the homogeneous heat equation such that
\begin{equation}\label{eq:Meier-global}
\int_0^{\infty} h(t)\,\|w(t,\cdot)\|_{L^\infty(\RD)}^{p-1}\,dt
\;<\;\infty,
\end{equation}
then the problem admits global positive solutions
{$u\le C w$}, with $\|u(t)\|_\infty\to 0$ as $t\to\infty$
{whenever $\|w(t)\|_\infty\to0$}.

\textup{(ii)} \emph{(Blow-up.)} If, for every nonnegative nontrivial
solution $w$ of the homogeneous heat equation,
\begin{equation}\label{eq:Meier-blow}
\limsup_{t\to\infty}
\|w(t,\cdot)\|_{L^\infty(\RD)}^{p-1}\int_0^{t}h(s)\,ds=\infty,
\end{equation}
then every nontrivial nonnegative solution blows up in finite time.
\end{theo}

\begin{theo}[Meier's universal formula {\cite[Theorem 1]{Meier90}}]
\label{thm:Meier1}
Assume \chk{hypotheses (h1)--(h2) of \cite{Meier90}}, that is,
$\alpha_0 t^q\le h(t)\le \alpha_1 t^q$ for $t$ large and some $q>-1$.
Then the critical exponent for the Cauchy problem
$\partial_t u-\Delta u=h(t) u^p$ is
$$
p^{\,\ast}\;=\;1+\frac{q+1}{s^{\,\ast}},
$$
where $s^{\,\ast}$ is given by~\eqref{eq:Meier-formula}.
For $D=\RD$, $s^{\,\ast}=N/2$, hence $p^{\,\ast}=1+2(q+1)/N$.
\end{theo}

The proof of Theorem~\ref{thm:Meier2}(i) relies on an explicit
supersolution $\bar u(x,t)=a(t)w(x,t)$, where $a$ solves the ODE
$$
a'(t)\;=\;h(t)\,\|w(t,\cdot)\|_{L^\infty}^{p-1}\,a(t)^p,\qquad a(0)=a_0>0.
$$
{Indeed, since $w^{p-1}\le\|w\|_\infty^{p-1}$, we get
$\bar u_t-\Delta\bar u=a'w\ge h\,a^pw^p=h\,\bar u^{\,p}$.}
Separating variables, we obtain
$a(t)^{1-p}=a_0^{1-p}-(p-1)\int_0^{t}h(s)\|w(s)\|_\infty^{p-1}ds$.
{This quantity stays bounded away from zero for all $t>0$, so that
$a$ remains bounded,} provided $a_0$ is small enough compared with the
finite integral~\eqref{eq:Meier-global}.

The blow-up part is proved by a subsolution argument. One constructs
$\underline u(x,t)=z(v(t);w(x,t))$, where $v(t)=\int_0^t h(s)\,ds$ and
$z(\cdot\,;\xi)$ solves the autonomous ODE $z'=z^p$ with $z(0;\xi)=\xi$.
{Since $z$ is convex in $\xi$, a direct computation gives
$\underline u_t-\Delta\underline u=h\,\underline u^{\,p}-z_{\xi\xi}|\nabla w|^2\le h\,\underline u^{\,p}$.}

\begin{rem}\rm\label{rem:no-leibniz}
\chk{Comparison principles for mild solutions remain valid for the
operator $\CapD{\alpha}+\Lap$, because all the kernels involved are
nonnegative} (see e.g.\ \cite{Zacher09,KilbasBook,CortazarQuirosWolanski24}).
What fails for $\alpha<1$ is the \emph{product ansatz}. The Caputo
derivative does not satisfy a Leibniz rule, so the supersolution
$\bar u=a(t)w(x,t)$ is no longer compatible with the nonlocal time
operator. The proofs must therefore be modified in the fractional
setting. {Global existence is established through fixed-point
arguments for mild solutions, and blow-up through Kaplan-type
functionals ($\alpha=1$) and nonlinear iterations on suitable weighted
averages ($\alpha<1$); see Section~\ref{sec:methods}.} The fate of Meier's
universal formula under this modification is discussed in
\S\ref{subsec:meier-frac}. At $\alpha=1$ it persists verbatim with
$s^\ast=N/\theta$, while for $\alpha<1$ the memory kernel shifts the
numerator $q+1$ to $q+\alpha$.
\end{rem}

\section{Main results}
\label{sec:main}

We consider the Cauchy problem~\eqref{eq:main} with the source $h$ as in
Convention~\ref{conv:source} and nonnegative initial data
$u_0\in L^1(\RD)\cap L^\infty(\RD)$, $u_0\not\equiv 0$. In the
global-existence statements, smallness is always measured in
$\|u_0\|_{L^1}+\|u_0\|_{L^\infty}$, while the blow-up statements hold for
\emph{every} nontrivial nonnegative datum in this class. Under
Convention~\ref{conv:source} the source is continuous up to $t=0$, so a
{mild solution in the sense of Definition~\ref{def:mild}} exists
locally and is unique. It is called \emph{global} if it can be extended
to all of $[0,\infty)$ and is bounded on every $[0,T]$. Otherwise we say
that the solution \emph{blows up in finite time}, which means that there
exists $T^{\,\ast}<\infty$ such that
$\limsup_{t\uparrow T^{\,\ast}}\|u(t,\cdot)\|_{L^\infty(\RD)}=\infty$.

Our first main result settles the case of the classical time derivative
completely.

\begin{theo}
\label{thm:main}
Let $\alpha=1$, $\theta\in(0,2]$ and $p>1$, and let
$\ell:(0,\infty)\to(0,\infty)$ be continuous and slowly varying at
infinity. Set $h(t)=t^{\gamma}\ell(t)$ for large $t$, in the sense of
Convention~\ref{conv:source}.
\begin{enumerate}
\item[\textup{(i)}] \emph{(No Fujita regime: $\gamma\le-1$.)}
For every $p>1$ and all sufficiently small nontrivial nonnegative data
$u_0\in L^1\cap L^\infty$, the mild solution of~\eqref{eq:main} is
global and satisfies $\|u(t,\cdot)\|_{L^\infty(\RD)}\to 0$ as
$t\to\infty$. There is no Fujita threshold in this regime; formally,
$p_c(\gamma)=\infty$.
\item[\textup{(ii)}] \emph{(Fujita regime: $\gamma>-1$.)} The critical
exponent is
$$
p_c(\gamma) \;=\; 1+\frac{\theta(\gamma+1)}{N}.
$$
If $1<p<p_c(\gamma)$, every nontrivial nonnegative solution blows up in
finite time. If $p>p_c(\gamma)$, small global solutions exist and satisfy
$\|u(t,\cdot)\|_{L^\infty(\RD)}\le C(1+t)^{-N/\theta}$.
\item[\textup{(iii)}] \emph{(Critical case: $\gamma>-1$, $p=p_c(\gamma)$.)}
The behavior at the critical exponent is determined by the slowly
varying function:
$$
\int_{1}^{\infty} t^{-1}\ell(t)\,dt =\infty
\;\;\Longrightarrow\;\; \text{every nontrivial solution blows up in finite time},
$$
$$
\int_{1}^{\infty} t^{-1}\ell(t)\,dt <\infty
\;\;\Longrightarrow\;\; \text{small global solutions exist}.
$$
\end{enumerate}
\end{theo}

For the genuinely fractional time derivative the exponents change. For
$\alpha\in(0,1)$ and $\gamma>-\alpha$ we set
\begin{equation}\label{eq:pcfrac}
p_c(\gamma,\alpha)\;:=\;1+\frac{\theta(\gamma+\alpha)}{N\alpha},
\end{equation}
and note that $p_c(\gamma,1)=p_c(\gamma)$ and $p_c(0,\alpha)=1+\theta/N$.

\begin{theo}
\label{thm:frac}
Let $\alpha\in(0,1)$, $\theta\in(0,2]$ and $p>1$, and let $h$ be as in
Convention~\ref{conv:source}.
\begin{enumerate}
\item[\textup{(a)}] \emph{(No Fujita regime: $\gamma\le-\alpha$.)}
For every $p>1$ and every $N\ge1$, small data in $L^1\cap L^\infty$
generate global mild solutions. In particular, global existence for all
$p>1$ extends from $\gamma\le-1$ to the strictly larger range
$\gamma\le-\alpha$.
\item[\textup{(b)}] \emph{(Global existence: $\gamma>-\alpha$.)}
Small-data global existence holds in the following cases.
\begin{itemize}
\item $N<\theta$ (since $N\ge 1$, this means $N=1$ and $\theta\in(1,2]$)
and $p>p_c(\gamma,\alpha)$. In this case the solutions satisfy
$\|u(t,\cdot)\|_{L^\infty(\RD)}\le C(1+t)^{-N\alpha/\theta}$.
Moreover, if $N<\theta$ and $\ell$ is bounded on $[1,\infty)$,
small-data global existence also holds at $p=p_c(\gamma,\alpha)$.
\item $N\ge\theta$ and $p>1+\dfrac{\gamma+\alpha}{\alpha}$.
\end{itemize}
Note that $1+(\gamma+\alpha)/\alpha\ge p_c(\gamma,\alpha)$ when
$N\ge\theta$, with equality if and only if $N=\theta$.
\item[\textup{(c)}] \emph{(Blow-up: $\gamma>-\alpha$, every $N\ge1$.)}
If $1<p<p_c(\gamma,\alpha)$, every nontrivial nonnegative solution blows
up in finite time. The same conclusion holds at $p=p_c(\gamma,\alpha)$
provided that $\limsup_{t\to\infty}\ell(t)=\infty$.
\end{enumerate}
\end{theo}

\begin{prop}
\label{prop:local}
Let $\alpha\in(0,1)$, $\theta\in(0,2]$, $N>\theta$, $\gamma>0$, and let
$h$ be as in Convention~\ref{conv:source}. Assume that either
$1<p<1+\gamma/\alpha$, or $p=1+\gamma/\alpha$ and
$\limsup_{t\to\infty}\ell(t)=\infty$. Then every nontrivial nonnegative
mild solution of~\eqref{eq:main} blows up in finite time.
\end{prop}

Combining parts (b) and (c) of Theorem~\ref{thm:frac}, we obtain a
complete picture in the one-dimensional regime $N<\theta$.

\begin{coro}\label{cor:N1}
Let $\alpha\in(0,1)$, $N=1$, $\theta\in(1,2]$ and $\gamma>-\alpha$. Then
$p_c(\gamma,\alpha)$ is the Fujita exponent of~\eqref{eq:main}. At
$p=p_c(\gamma,\alpha)$, small global solutions exist if $\ell$ is bounded
on $[1,\infty)$, and every nontrivial solution blows up if $\ell$ is
unbounded on $[1,\infty)$.
\end{coro}

Since $\ell$ is continuous, it is unbounded on $[1,\infty)$ if and only
if $\limsup_{t\to\infty}\ell(t)=\infty$, so the two alternatives in
Corollary~\ref{cor:N1} are exhaustive.

For $N>\theta$, Theorem~\ref{thm:frac} and Proposition~\ref{prop:local}
leave a gap between the blow-up range
$1<p<\max\{p_c(\gamma,\alpha),1+\gamma/\alpha\}$ and the global range
$p>1+(\gamma+\alpha)/\alpha$. Note that
$1+\gamma/\alpha>p_c(\gamma,\alpha)$ if and only if
$\gamma>\alpha\theta/(N-\theta)$; in this regime $p_c(\gamma,\alpha)$ is
\emph{not} the Fujita exponent.

\begin{conj}\label{conj:frac}
Let $\alpha\in(0,1)$, $N>\theta$ and $\gamma>-\alpha$. Then the Fujita
exponent of~\eqref{eq:main} is
$$
p_\ast(\gamma,\alpha)\;:=\;\max\Bigl\{p_c(\gamma,\alpha),\;1+\frac{\gamma}{\alpha}\Bigr\},
$$
that is, small global solutions exist for every $p>p_\ast(\gamma,\alpha)$.
For $N=\theta$, small global solutions exist at $p=p_c(\gamma,\alpha)$
whenever $\ell$ is bounded.
\end{conj}
The first term in $p_\ast$ comes from the scaling (mass) mechanism of
Theorem~\ref{thm:frac}(c), and the second one from the local mechanism of
Proposition~\ref{prop:local}. The conjecture is supported only by
heuristics: on balls of fixed size, the time-fractional dynamics behaves
like the damped fractional ODE
$\CapD{\alpha}v+\lambda v=t^{\gamma}v^p$, whose small solutions decay like
$t^{-\alpha}$ and remain global precisely when $p>1+\gamma/\alpha$. For
$\gamma=0$ it reduces to the results of
\cite{Zhang14,CortazarQuirosWolanski24}.

A few remarks are in order.
\begin{rem}\rm
\label{rem:shift-proved}
Meier's classification cannot persist verbatim for $\alpha\in(0,1)$, and
Theorem~\ref{thm:frac}(a) alone already shows it. Fix $\alpha\in(0,1)$
and $\gamma\in(-1,-\alpha]$, which is a non-empty range. Any
classification driven by the numerator $\gamma+1$, that is, any
candidate critical exponent of the form $1+c\,(\gamma+1)$ with $c>0$ as
produced by Meier's formula~\eqref{eq:Meier-formula} for \emph{any}
choice of $s^{\,\ast}$, would be finite and larger than $1$ for such
$\gamma$. It would therefore predict a non-empty blow-up window
$1<p<1+c\,(\gamma+1)$. Theorem~\ref{thm:frac}(a) excludes this: for every
$\gamma\le-\alpha$, small-data global solutions exist for \emph{every}
$p>1$, although the source $t^{\gamma}\ell(t)$ is not integrable at
infinity when $\gamma\in(-1,-\alpha]$. Hence, for $\alpha\in(0,1)$, the
Fujita regime starts at $\gamma=-\alpha$ and not at $\gamma=-1$.
{Conversely, the blow-up range $1<p<p_c(\gamma,\alpha)$ of
Theorem~\ref{thm:frac}(c) shows that the shift also occurs on the
blow-up side, for every $\gamma>-\alpha$ and every $N\ge1$.} The memory of
the Caputo derivative thus \emph{provably} shifts the numerator from
$\gamma+1$ to $\gamma+\alpha$.
\end{rem}

\begin{rem}\rm
When $\theta=2$ (and $\alpha=1$), Theorem~\ref{thm:main}(i)--(ii)
recovers the results of Meier~\cite{Meier90} on $\RD$, with
$p_c(\gamma)=1+2(\gamma+1)/N$.
{The critical-case dichotomy (iii) can be read off the integral
$\int^{\infty} h(t)\,t^{-N(p-1)/\theta}\,dt$, which, for
$h(t)=t^{\gamma}\ell(t)$ and $p=p_c(\gamma)$, reduces exactly to
$\int^\infty t^{-1}\ell(t)\,dt$. Its divergence is the blow-up criterion of
Guedda and Kirane~\cite{GueddaKirane}, while its convergence is precisely
condition~\eqref{eq:Meier-global} for $w=S_{1,\theta}(t)w_0$, because
$\|S_{1,\theta}(t)w_0\|_\infty\asymp t^{-N/\theta}$ for large $t$.
Meier's supersolution $a(t)w$ is also a supersolution when $-\Delta$ is
replaced by $\Lap$. Consequently, at $\alpha=1$ the individual ingredients
of Theorem~\ref{thm:main} are essentially available in the
literature~\cite{Meier90,GueddaKirane,Sugitani75}. Our contribution there
is a unified treatment of all $\theta\in(0,2]$ and the formulation of the
critical case in terms of Karamata's theory. The genuinely new results
concern $\alpha\in(0,1)$, namely Theorem~\ref{thm:frac},
Proposition~\ref{prop:local} and Corollary~\ref{cor:N1}.}
\end{rem}

\begin{rem}\rm
The dependence of the critical exponent on the time order $\alpha$ is
more subtle than the autonomous case suggests. Writing
$p_c(\gamma,\alpha)-1=\frac{\theta}{N}\bigl(1+\frac{\gamma}{\alpha}\bigr)$,
we see the following. At $\gamma=0$, the exponent $1+\theta/N$ does not
depend on $\alpha$, in agreement with
\cite{Zhang14,CortazarQuirosWolanski24}. For $\gamma>0$, a smaller
$\alpha$ (more memory) yields a \emph{larger} exponent, and hence a larger
blow-up window. For $\gamma\in(-\alpha,0)$ the monotonicity is reversed.
Finally, the onset of the Fujita regime itself moves from $\gamma=-1$
(at $\alpha=1$) to $\gamma=-\alpha$. Heuristically, the memory kernel
$(t-s)^{\alpha-1}$ in~\eqref{eq:mild} converts the cumulative source
$\int_0^t h(s)\,ds$ of Meier's criterion into the weighted average
$\int_0^t (t-s)^{\alpha-1}h(s)\,ds\sim B(\alpha,\gamma+1)\,t^{\gamma+\alpha}\ell(t)$
{(valid for $\gamma>-1$)}, which replaces $\gamma+1$ by $\gamma+\alpha$
throughout. {For $N>\theta$, the local mechanism of
Proposition~\ref{prop:local} adds a second, $\theta$-independent
constraint $p>1+\gamma/\alpha$.}
\end{rem}

\begin{rem}\rm
The slowly varying factor $\ell$ has \emph{no effect} on the value of the
critical exponent. At $\alpha=1$, it is decisive exactly \emph{at} the
critical exponent $p=p_c(\gamma)$, for every $\gamma>-1$. There, the
dichotomy is controlled by the single Karamata-type integral
$\int_1^\infty t^{-1}\ell(t)\,dt$, and Lemma~\ref{lem:karamata}(K3)
yields the following examples.
\begin{itemize}
\item $\ell\equiv 1$: divergence, hence blow-up at criticality (for
$\gamma=0$, $\theta=2$ this is the classical critical-case blow-up of
Hayakawa~\cite{Hayakawa73} and Kobayashi, Sirao and Tanaka~\cite{KST77});
\item $\ell(t)=(\log t)^{a}$ for $t\ge e$: divergence if and only if
$a\ge-1$, convergence if and only if $a<-1$;
\item $\ell(t)=(\log t)^{-1}(\log\log t)^{b}$ for large $t$: divergence
if and only if $b\ge-1$.
\end{itemize}
{For $\alpha\in(0,1)$ and $N<\theta$, Corollary~\ref{cor:N1} shows
that the relevant property is the boundedness of $\ell$, a weaker
requirement than the convergence of the Karamata integral. For instance,
for $\ell(t)=(\log t)^{-1/2}$ (for $t\ge e$), solutions blow up at the
critical exponent when $\alpha=1$, whereas small global solutions exist at
the critical exponent when $\alpha<1$. The memory of the Caputo
derivative thus turns critical blow-up into critical global existence for
every bounded $\ell$, in agreement with the autonomous case $\ell\equiv1$,
$\gamma=0$ \cite{Zhang14,CortazarQuirosWolanski24}.}
\end{rem}

\begin{rem}\rm
The dichotomy in Theorem~\ref{thm:main}(iii) is sharp. At
$p=p_c(\gamma)$, \emph{any} slowly varying $\ell$ with a finite tail
integral $\int_1^{\infty}t^{-1}\ell(t)\,dt$ allows small global
solutions, while any $\ell$ with a divergent tail integral forces blow-up
of every nontrivial solution. In particular, an arbitrarily mild
logarithmic damping, for instance $\ell(t)=(\log t)^{-1-\delta}$ with
$\delta>0$, suffices to move the critical exponent to the
global-existence side.
\end{rem}

\begin{example}\rm
Let $\alpha=1$, fix any $\gamma>-1$, and take $\ell(t)=(\log t)^{-2}$
for $t\ge e$ {(extended to $(0,e)$ as a positive continuous
function)}. Then
$$
\int_{e}^{\infty} \frac{1}{t\,(\log t)^{2}}\,dt
\;=\;\Bigl[-\frac{1}{\log t}\Bigr]_{e}^{\infty}\;=\;1\;<\;\infty.
$$
By Theorem~\ref{thm:main}(iii), small global solutions exist \emph{at}
the critical exponent $p=p_c(\gamma)=1+\theta(\gamma+1)/N$. This is
impossible for $\ell\equiv 1$, where the critical case always belongs to
the blow-up regime \cite{Hayakawa73,KST77,Sugitani75}.
\end{example}

\begin{example}\rm
Let $\alpha=1$, $\gamma=0$ and $\ell(t)=(\log(2+t))^{a}$ with $a\in\R$.
This corresponds to $h(t)\sim(\log t)^a$, a logarithmic perturbation of
the autonomous source. Theorem~\ref{thm:main}(ii) gives
$p_c(0)=1+\theta/N$ independently of $a$, so the logarithmic correction
does \emph{not} shift the threshold. It does, however, decide the
critical case. At $p=1+\theta/N$, every nontrivial solution blows up when
$a\ge -1$, whereas small global solutions exist when $a<-1$.
\end{example}

\begin{example}\rm
Let $\alpha=1$ and $\gamma>-1$, and consider the doubly logarithmic
correction
$\ell(t)=\bigl(\log(2+t)\bigr)^{a}\bigl(\log\log(16+t)\bigr)^{b}$ with
$a,b\in\R$. The critical exponent is still
$p_c(\gamma)=1+\theta(\gamma+1)/N$, independently of $(a,b)$. At
$p=p_c(\gamma)$, blow-up occurs precisely when $a>-1$, or when $a=-1$
and $b\ge-1$.
\end{example}

\section{Other classes of time-dependent coefficients}
\label{sec:other}

We briefly comment on classes of source intensities $h$ that are not
covered by the Karamata framework $h(t)=t^\gamma\ell(t)$.

\subsection{Exponential growth}
\label{subsec:exp}

Let $h(t)=e^{\mu t}$ with $\mu>0$; the growth rate $\mu$ is unrelated
to the decay exponent $\beta=N\alpha/\theta$ of
\S\ref{subsec:proof-existence}. The ODE $a'(t)=e^{\mu t}a(t)^p$, which is
the analogue of Meier's auxiliary equation, has the explicit solution
$$
a(t)^{1-p}\;=\;a_0^{1-p}-\frac{p-1}{\mu}\bigl(e^{\mu t}-1\bigr),
$$
which blows up at
$$
T^{\,\ast}\;=\;\mu^{-1}\log\!\left(1+\frac{\mu\,a_0^{1-p}}{p-1}\right)
\;=\;\mu^{-1}\log\!\left(1+\frac{\mu}{(p-1)\,a_0^{p-1}}\right).
$$
Thus $T^{\,\ast}$ is finite for every $a_0>0$, and $T^{\,\ast}\to\infty$
as $a_0\to 0^+$. Hence, even without diffusion, the source alone forces
blow-up in finite time for every $p>1$. In the PDE setting on the whole
space $\RD$ with $\alpha=1$, the extreme growth of $h$ overwhelms the
dissipation. The Kaplan-type functional of \S\ref{subsec:proof-blowup}
yields the differential inequality $\psi'\ge e^{\mu t}\psi^p$, and the
resulting necessary condition for global existence,
$\sup_{T}T^{-N(p-1)/\theta}e^{\mu T}<\infty$, fails for every $p>1$.
Consequently, \emph{every nontrivial nonnegative solution blows up in
finite time for every $p>1$}, regardless of $\theta$. The same behavior
is expected for $\alpha<1$, although the memory term precludes the simple
adjoint argument. On a bounded domain the situation is more interesting;
see Meier~\cite[Theorem 4]{Meier90} for the sharp threshold
$p^{\,\ast}_\mu=1+\mu/\lambda_1$, where $\lambda_1$ is the first
Dirichlet eigenvalue of the Laplacian.

\subsection{Pure power behavior}
\label{subsec:power}

If $h(t)=t^q$, with no slowly varying correction, we are in the regularly
varying case with $\ell\equiv 1$. For $\alpha=1$ and $q>-1$,
Theorem~\ref{thm:main} yields $p_c=1+\theta(q+1)/N$, with blow-up
\emph{including} the critical exponent, since
$\int^{\infty}t^{-1}\,dt=\infty$. Taking $\theta=2$ recovers Meier's
value $p_c=1+2(q+1)/N$ \cite{Meier90}, while $\theta\in(0,2)$ gives its
natural $\theta$-stable analogue, which reduces to Sugitani's
theorem~\cite{Sugitani75} at $q=0$. For $\alpha<1$, see
Theorem~\ref{thm:frac}, {Proposition~\ref{prop:local}} and
Conjecture~\ref{conj:frac}.

\subsection{Oscillatory coefficients}
\label{subsec:osc}

Slowly varying functions may exhibit persistent oscillations. A
canonical example is $\ell(t)=2+\sin(\log\log t)$ for large $t$. It is
slowly varying because
$$
\log\log(\lambda t)-\log\log t
\;=\;\log\!\Bigl(1+\frac{\log\lambda}{\log t}\Bigr)\;\longrightarrow\; 0
\qquad(t\to\infty),
$$
so that $\sin(\log\log(\lambda t))-\sin(\log\log t)\to 0$ for every
$\lambda>0$, and yet $\ell$ does not converge as $t\to\infty$. Note that
the often-quoted candidate $2+\sin(\log t)$ is \emph{not} slowly varying.
Indeed,
$\sin(\log(\lambda t))-\sin(\log t)
=2\cos\!\bigl(\log t + \tfrac{\log\lambda}{2}\bigr)\sin\!\bigl(\tfrac{\log\lambda}{2}\bigr)$
oscillates between $\pm 2|\sin(\log\lambda/2)|$ as $t\to\infty$, for
every fixed $\lambda\ne 1$ with $\sin(\log\lambda/2)\neq0$.
Theorem~\ref{thm:main} still applies. The value of $p_c(\gamma)$ does
not see such oscillations, and at $p=p_c(\gamma)$ the dichotomy depends
only on the Karamata integral $\int_1^{\infty}t^{-1}\ell(t)\,dt$. For
$\ell(t)=2+\sin(\log\log t)$ this integral diverges (since $\ell\ge 1$),
so blow-up persists at criticality. {For $\alpha<1$ and $N<\theta$,
this $\ell$ is bounded, and Corollary~\ref{cor:N1} gives global existence
at the critical exponent.}

\subsection{Integrable sources of non-power decay}
\label{subsec:int}

{Assume that $\alpha=1$, that $h\in C([0,\infty))$ is nonnegative and
$h\in L^1(0,\infty)$, but that $h$ is \emph{not} regularly varying; for
example $h(t)=e^{-t}$ or $h(t)=(1+t)^{-2}(2+\sin t)$. Integrability alone
then forces global existence for small data, for every $p>1$. Indeed, in
the Duhamel fixed-point argument of \S\ref{subsec:proof-existence}, the
contribution of $[0,t/2]$ is bounded by
$C\,t^{-\beta}\int_0^\infty h$, and that of $[t/2,t]$ by
$C(1+t)^{-\beta p}\int_0^\infty h$. For $\alpha<1$, integrability alone is
not sufficient for this argument, since the near-diagonal term
$\int_{t/2}^t(t-s)^{\alpha-1}h(s)\,ds$ need not decay. A pointwise bound
such as $h(t)\le C(1+t)^{-\alpha}$ restores the conclusion through
Theorem~\ref{thm:frac}(a) and its proof.} The autonomous-case analysis is
then irrelevant.

\subsection{More general regularly varying functions}
\label{subsec:gen}

If $h(t)=t^\gamma g(t)$ with $g$ \emph{regularly varying} of index
$\rho\ne 0$, i.e.\ $g(\lambda t)/g(t)\to\lambda^\rho$ for every
$\lambda>0$, then $h$ is itself regularly varying of index $\gamma+\rho$.
At $\alpha=1$, the critical exponent formula becomes
$$
p_c(\gamma+\rho)\;=\;1+\frac{\theta(\gamma+\rho+1)}{N},
$$
provided that $\gamma+\rho>-1$. The critical-case dichotomy is then
governed by the slowly varying function $g(t)t^{-\rho}$, through the
integral $\int^\infty t^{-1}\,g(t)t^{-\rho}\,dt$.

\subsection{Two-sided power bounds (Meier's assumption (h2))}
\label{subsec:h2}

Meier~\cite{Meier90} works under the weaker assumption
$$
\alpha_0\, t^{q}\;\le\; h(t)\;\le\;\alpha_1\, t^{q}\qquad\text{for }t\text{ large},
$$
with constants $0<\alpha_0\le\alpha_1$ (here $\alpha_0,\alpha_1$ are
unrelated to the Caputo order $\alpha$). Our framework covers this case
with $\gamma=q$, since such an $h$ is sandwiched between two regularly
varying functions with the same index. The global-existence estimates of
\S\ref{subsec:proof-existence} use only the upper bound, and the blow-up
arguments of \S\S\ref{subsec:proof-blowup}--\ref{subsec:local} use only
the lower bound. Hence the proofs adapt directly away from the critical
exponent. At the critical exponent, the two-sided bound does not
determine the dichotomy, which is genuinely a property of~$\ell$.

\section{Proofs of the main results}
\label{sec:methods}

{This section contains the proofs. Global existence is obtained by a
fixed-point argument in weighted spaces
(\S\ref{subsec:proof-existence}). For blow-up at $\alpha=1$ we use a
Kaplan-type functional built on the kernel $Z$, together with the
integral criterion of Guedda and Kirane~\cite{GueddaKirane} at the
critical exponent (\S\ref{subsec:proof-blowup}). For $\alpha\in(0,1)$ we
adapt the nonlinear iteration of Cort\'azar, Quir\'os and
Wolanski~\cite{CortazarQuirosWolanski24} to the time-dependent source
(\S\ref{subsec:proof-blowup-frac}), and we prove the local blow-up
mechanism of Proposition~\ref{prop:local} in \S\ref{subsec:local}. The
regime $\gamma\le-\alpha$ and the fate of Meier's formula are discussed in
\S\S\ref{subsec:noFujita}--\ref{subsec:meier-frac}. Throughout, the
kernels and solutions are written as functions of $(t,x)$, in this order.}

\subsection{Global existence}
\label{subsec:proof-existence}

Throughout this subsection we work under the standing hypothesis
\begin{equation}\label{eq:hypS}
\textup{(S)}\qquad \alpha=1\ \text{(any $N\ge1$)}
\qquad\text{or}\qquad \alpha\in(0,1)\ \text{and}\ N<\theta.
\end{equation}
By Proposition~\ref{prop:heat-kernel} and
Remark~\ref{rem:kernel-singularity}, this is exactly the regime in which
the $L^1\to L^\infty$ smoothing estimate
$\|S_{\alpha,\theta}(t)\|_{L^1\to L^\infty}\le Ct^{-N\alpha/\theta}$
holds. {Under (S), the analogous estimate
$\|P_{\alpha,\theta}(t)\|_{L^1\to L^\infty}\le Ct^{-N\alpha/\theta}$ also
holds, because $N<\theta<2\theta$ when $\alpha<1$ (and
$P_{1,\theta}=S_{1,\theta}$).} The case $\alpha<1$, $N\ge\theta$ is
treated in Remark~\ref{rem:Lr-scheme} below.

\medskip
\noindent\textbf{Set-up.} Set $\beta:=N\alpha/\theta$ and let
{
$$
X_{\beta}\;:=\;\bigl\{u:(0,\infty)\times\RD\to\R \text{ measurable}\,:\,
\|u\|_{X_\beta}<\infty\bigr\},
$$
endowed with the norm
$$
\|u\|_{X_\beta}\;:=\;\operatorname*{ess\,sup}_{t>0}\Bigl[\|u(t,\cdot)\|_{L^1(\RD)}\;+\;(1+t)^{\beta}\,\|u(t,\cdot)\|_{L^\infty(\RD)}\Bigr].
$$
Then $X_\beta$ is a Banach space, and its nonnegative elements form a
closed subset $X_\beta^+$.} For nonnegative
$u_0\in L^1(\RD)\cap L^\infty(\RD)$, we define the map
$$
\Lambda[u](t)\;=\;S_{\alpha,\theta}(t)u_0
\;+\;\int_0^{t}(t-s)^{\alpha-1}P_{\alpha,\theta}(t-s)\bigl[h(s)\,u(s)^p\bigr]\,ds=:\mathcal{L}[u](t)+\mathcal{N}[u](t).
$$
A mild solution of~\eqref{eq:main} is a fixed point of $\Lambda$.
{Since the kernels of $S_{\alpha,\theta}$ and $P_{\alpha,\theta}$ and
the function $h$ are nonnegative, $\Lambda$ maps $X_\beta^+$ into
nonnegative functions.}

\medskip
\noindent\textbf{Decay of the linear part.} Under \textup{(S)},
Proposition~\ref{prop:heat-kernel} with $(p,q)=(1,\infty)$ gives
$\|\mathcal{L}[u](t)\|_\infty\le C\, t^{-N\alpha/\theta}\|u_0\|_1$ for
$t\ge 1$, while $\|\mathcal{L}[u](t)\|_\infty\le\|u_0\|_\infty$ for
$t\in(0,1]$. Hence
$$\|\mathcal{L}[u](t)\|_\infty\le C_0(1+t)^{-N\alpha/\theta}\bigl(\|u_0\|_1+\|u_0\|_\infty\bigr).$$
The mass bound $\|\mathcal{L}[u](t)\|_1\le \|u_0\|_1$ follows at once
from the $L^1$ contractivity of $S_{\alpha,\theta}(t)$.

\medskip
\noindent\textbf{Nonlinear estimates.} For $u\in X_\beta^+$, we write
$M:=\|u\|_{X_\beta}$. Then, for almost every $s>0$,
\begin{equation}
\|(u(s))^p\|_{L^1(\RD)}\le \|u(s)\|_{L^\infty}^{p-1}\,\|u(s)\|_{L^1}\;\le\;M^{p}\,(1+s)^{-\beta(p-1)},
\end{equation}
and
\begin{equation}\label{estiC}
\|(u(s))^p\|_{L^\infty(\RD)}\;\le\;\|u(s)\|_{L^\infty}^p\;\le\;M^{p}\,(1+s)^{-\beta p}.
\end{equation}
The decisive quantity is
\begin{equation}\label{eq:kappa}
\kappa\;:=\;\gamma+\alpha-\beta(p-1),
\end{equation}
which measures the competition between the cumulative source, weighted
by the memory kernel, and the dispersive decay. With
$\beta=N\alpha/\theta$, we have
$$
\kappa<0 \iff p>1+\frac{\gamma+\alpha}{\beta}
=1+\frac{\theta(\gamma+\alpha)}{N\alpha}.
$$
At $\alpha=1$ this is exactly $p>p_c(\gamma)$, and it holds for
\emph{every} $p>1$ when $\gamma\le-\alpha$. For $\alpha<1$ and
$\gamma>-\alpha$ it reads $p>p_c(\gamma,\alpha)$.

\medskip
\noindent\textbf{The strictly supercritical regime $\kappa<0$.}
Fix $\varepsilon\in(0,-\kappa)$. By Lemma~\ref{lem:karamata}(K1) there
exists $t_0\ge 2$ such that
\begin{equation}\label{estiell}
   \ell(s)\leq C\,s^{\varepsilon},\qquad \hbox{for all}\,\, s\ge t_0.
\end{equation}

$\bullet$ {\bf Case $t\leq 2t_0$.} Using the uniform bound
$\|P_{\alpha,\theta}(t-s)f\|_\infty\le C\|f\|_\infty$, we get
\begin{eqnarray*}
\Bigl\|\mathcal{N}[u](t)\Bigr\|_\infty
&\le& C\,M^{p}\!\int_0^{t}(t-s)^{\alpha-1}\,h(s)\,(1+s)^{-\beta p}\,ds\\
&\leq &C\,M^{p}\!\int_0^{t}(t-s)^{\alpha-1}\,h(s)\,ds\\
&\le&C\,M^{p}(2t_0)^{\alpha},
\end{eqnarray*}
where we used the boundedness of $h$ on $[0,2t_0]$. Since
$(1+t)^{-\beta}\ge (1+2t_0)^{-\beta}$ for $t\le2t_0$, it follows that
\begin{equation}\label{eq:estiA}
\Bigl\|\mathcal{N}[u](t)\Bigr\|_\infty\le C\,M^{p}(1+t)^{-\beta}.
\end{equation}

$\bullet$ {\bf Case $t\geq 2t_0$.} We split the integral defining
$\mathcal{N}[u](t)$ into the intervals $[0,t/2]$ and $[t/2,t]$. On
$[0,t/2]$ we apply the $L^1$--$L^\infty$ estimate for $P_{\alpha,\theta}$
(cf.~\cite[Theorem~3.2]{CortazarQuirosWolanski24}),
$$\|P_{\alpha,\theta}(t-s)f\|_\infty\le C(t-s)^{-N\alpha/\theta}\|f\|_1,$$
which is valid under assumption \textup{(S)}. On $[t/2,t]$ we use the
uniform $L^\infty$ bound. We obtain
\begin{eqnarray*}
\Bigl\|\mathcal{N}[u](t)\Bigr\|_\infty
&\le& C\!\int_0^{t/2}(t-s)^{\alpha-1-\beta}\,h(s)\|(u(s))^p\|_{L^1}\,ds+C\!\int^t_{t/2}(t-s)^{\alpha-1}\,h(s)\|(u(s))^p\|_{\infty}\,ds\\
&\le& C\,M^p\bigl(I_1(t)+I_2(t)\bigr),
\end{eqnarray*}
where
$$I_1(t):=\int_0^{t/2}(t-s)^{\alpha-1-\beta}\,h(s)\,(1+s)^{-\beta(p-1)}\,ds
\quad\text{and}\quad
I_2(t):=\int^t_{t/2}(t-s)^{\alpha-1}\,h(s)\,(1+s)^{-\beta p}\,ds.$$
We first estimate $I_1(t)$. Since $t-s\ge t/2$ on $[0,t/2]$ and
$\alpha-1-\beta<0$, we have
$$I_1(t)\le C\,t^{\alpha-1-\beta}\int_0^{t/2}\,h(s)\,(1+s)^{-\beta(p-1)}\,ds.$$
We split the last integral at $t_0$. As $h$ is bounded on $[0,t_0]$ and
$t_0$ is fixed,
$$\int_0^{t_0}\,h(s)\,(1+s)^{-\beta(p-1)}\,ds\leq C(t_0)\;=:\;C.$$
On the other hand, $h(s)=s^{\gamma}\ell(s)$ for $s\ge t_0$, and
$\gamma+\varepsilon-\beta(p-1)=\kappa+\varepsilon-\alpha<-\alpha$.
{Using~\eqref{estiell}, we infer that
$$
\int_{t_0}^{t/2}\,h(s)\,(1+s)^{-\beta(p-1)}\,ds
\;\leq\; C \int_{t_0}^{t/2}\,s^{\gamma+\varepsilon-\beta(p-1)}\,ds
\;\leq\; C\,\bigl(1+t^{1-\alpha}\bigr).
$$
Indeed, set $\sigma:=\gamma+\varepsilon-\beta(p-1)<-\alpha$. If
$\sigma<-1$, the integral is bounded. If $\sigma=-1$, which can happen
only when $\alpha<1$, it is bounded by $\log t\le C\,t^{1-\alpha}$. If
$-1<\sigma<-\alpha$, it is bounded by $C\,t^{\sigma+1}\le C\,t^{1-\alpha}$.}
For $\alpha=1$ the integral is in fact convergent and the bound is a
constant. In either case
$$I_1(t)\leq C\,t^{\alpha-1-\beta}\bigl(1+t^{1-\alpha}\bigr)\leq C\, t^{-\beta}\leq C\,(1+t)^{-\beta},$$
where we used $\alpha\le 1$ and $t\ge 1$. Next, we estimate $I_2(t)$,
again using $h(s)=s^{\gamma}\ell(s)$ and~\eqref{estiell} for
$s\ge t/2\geq t_0$, together with $s\asymp t$ on $[t/2,t]$:
\begin{eqnarray*}
I_2(t)&=&\int_{t/2}^t(t-s)^{\alpha-1}s^{\gamma}\ell(s)\,(1+s)^{-\beta p}\,ds\\
&\leq& C \int_{t/2}^t(t-s)^{\alpha-1}\,s^{\gamma+\varepsilon-\beta p}\,ds\\
&\leq& C\,t^{\gamma+\varepsilon-\beta p} \int_{t/2}^t(t-s)^{\alpha-1}\,ds\\
&\le& C\,t^{\gamma+\varepsilon-\beta p+\alpha}\;=\;C\,t^{\kappa+\varepsilon}\,t^{-\beta}\\
&\leq& C\, t^{-\beta}\;\leq\; C\,(1+t)^{-\beta},
\end{eqnarray*}
since $\kappa+\varepsilon<0$. Combining the estimates of $I_1(t)$ and
$I_2(t)$ with~\eqref{eq:estiA}, we conclude that, for all $t>0$,
\begin{equation}
\Bigl\|\mathcal{N}[u](t)\Bigr\|_\infty\le C\,M^{p}(1+t)^{-\beta}.
\end{equation}
To estimate $\bigl\|\mathcal{N}[u](t)\bigr\|_{L^1}$, we use the $L^1$
bound $\|P_{\alpha,\theta}(t-s)f\|_1\le C\|f\|_1$ and get
$$
\Bigl\|\mathcal{N}[u](t)\Bigr\|_1
\le C\,M^{p}\!\int_0^{t}(t-s)^{\alpha-1}\,h(s)\,(1+s)^{-\beta(p-1)}\,ds.
$$
{For $t\le 2t_0$ the right-hand side is bounded by $C M^p$, since $h$
is bounded on $[0,2t_0]$.} For $t\ge2t_0$ we split at $t/2$ exactly as
above. The contribution of $[0,t/2]$ is at most
$C\,t^{\alpha-1}(1+t^{1-\alpha})\le C$, while that of $[t/2,t]$ is at most
$C\,t^{\gamma+\varepsilon-\beta(p-1)+\alpha}=C\,t^{\kappa+\varepsilon}\le C$.
Hence
\begin{equation}
\Bigl\|\mathcal{N}[u](t)\Bigr\|_{L^1}\le C\,M^{p}.
\end{equation}

\medskip
\noindent\textbf{The marginal regime $\kappa=0$.} Here there is no room
left for $\varepsilon$, and the slowly varying factor itself decides
whether the scheme closes. We treat the classical and the genuinely
fractional cases separately.

$\bullet$ \emph{$\alpha=1$, $p=p_c(\gamma)$, under
$\int_1^{\infty}t^{-1}\ell(t)\,dt<\infty$.} Since $\beta(p-1)=\gamma+1>0$,
the large-time density satisfies
$h(s)(1+s)^{-\beta(p-1)}\le s^{-1}\ell(s)$, so that
$$
\int_{t_0}^{t/2}h(s)(1+s)^{-\beta(p-1)}\,ds
\;\le\;\int_{t_0}^{\infty}s^{-1}\ell(s)\,ds\;<\;\infty,
$$
and $I_1(t)\le C\,t^{\alpha-1-\beta}=C\,t^{-\beta}$. For $I_2$ we cannot
afford any positive power. Instead, by Lemma~\ref{lem:karamata}(K4), the
convergence of the borderline integral forces $\ell(t)\to 0$. Hence
$\ell$ is bounded on $[t_0,\infty)$ and
$$
I_2(t)\;\le\;C\Bigl(\sup_{[t/2,t]}\ell\Bigr)\,t^{\gamma-\beta p}\,t^{\alpha}
\;=\;C\Bigl(\sup_{[t/2,t]}\ell\Bigr)\,t^{\kappa}\,t^{-\beta}\;\le\;C\,t^{-\beta}.
$$
The $L^1$ estimate is handled in the same way. Thus the scheme closes
\emph{at} the critical exponent, which proves the global-existence half of
Theorem~\ref{thm:main}(iii).

$\bullet$ \emph{$\alpha\in(0,1)$, $N<\theta$, $p=p_c(\gamma,\alpha)$,
with $\ell$ bounded.} Now $\beta(p-1)=\gamma+\alpha$, and the large-time
density is bounded by $s^{-\alpha}\ell(s)\le C s^{-\alpha}$. Therefore
$\int_{t_0}^{t/2}s^{-\alpha}\ell(s)\,ds\le C\,t^{1-\alpha}$ and
$I_1(t)\le C\,t^{-\beta}$, while
$I_2(t)\le C(\sup\ell)\,t^{\kappa}t^{-\beta}=C\,t^{-\beta}$. The $L^1$
estimate is again identical. This proves the critical-case global
existence in Theorem~\ref{thm:frac}(b), which is consistent with the
autonomous critical-case global existence of
\cite{Zhang14,CortazarQuirosWolanski24} ($\gamma=0$, $\ell\equiv1$).

\medskip
\noindent\textbf{Fixed-point argument.}
{The estimates above show that there is a constant $C_\ast$,
independent of $u_0$, such that
$\|\mathcal{N}[u]\|_{X_\beta}\le C_\ast\|u\|_{X_\beta}^p$ for every
$u\in X_\beta^+$. Moreover, the elementary inequality
$|a^p-b^p|\le p\,(a\vee b)^{p-1}|a-b|$ for $a,b\ge0$ gives
$\|u(s)^p-v(s)^p\|_{L^1}\le pM^{p-1}(1+s)^{-\beta(p-1)}\|u-v\|_{X_\beta}$
and
$\|u(s)^p-v(s)^p\|_{L^\infty}\le pM^{p-1}(1+s)^{-\beta p}\|u-v\|_{X_\beta}$
whenever $\|u\|_{X_\beta},\|v\|_{X_\beta}\le M$. These are the same
densities as above, so that
$\|\mathcal{N}[u]-\mathcal{N}[v]\|_{X_\beta}\le pC_\ast M^{p-1}\|u-v\|_{X_\beta}$.
Let $M:=2C_0(\|u_0\|_1+\|u_0\|_\infty)$, and assume that $u_0$ is so small
that $pC_\ast M^{p-1}\le 1/2$. Then, for $u$ in the closed set
$B_M^+:=\{u\in X_\beta^+:\|u\|_{X_\beta}\le M\}$, we have
$\|\Lambda[u]\|_{X_\beta}\le M/2+C_\ast M^{p}\le M$. Thus $\Lambda$ maps
$B_M^+$ into itself and is a contraction there, with constant $1/2$.} The
Banach fixed-point theorem produces a unique fixed point. By uniqueness of
mild solutions, it coincides with the maximal mild solution, which is
therefore global and satisfies
$\|u(t)\|_\infty\le M(1+t)^{-\beta}\to 0$ as $t\to\infty$.

\medskip
\noindent\textbf{Summary of the closure region.} The scheme closes, and
hence small-data global existence holds with decay $(1+t)^{-N\alpha/\theta}$
under hypothesis~\textup{(S)}, precisely when
$$
\kappa<0,
\qquad\text{or}\qquad
\kappa=0\ \text{ with }\
\begin{cases}
\int^{\infty}t^{-1}\ell(t)\,dt<\infty, & \alpha=1,\\[2pt]
\ell\ \text{bounded}, & \alpha<1.
\end{cases}
$$
At $\alpha=1$ this gives exactly all $p>1$ for $\gamma\le-1$,
$p>p_c(\gamma)$ for $\gamma>-1$, and $p=p_c(\gamma)$ under the convergent
Karamata integral; these are the global halves of
Theorem~\ref{thm:main}. For $\alpha<1$ (with $N<\theta$) it gives all
$p>1$ for $\gamma\le-\alpha$, and $p>p_c(\gamma,\alpha)$ together with
the marginal case with bounded $\ell$ for $\gamma>-\alpha$; these are
Theorem~\ref{thm:frac}(a)--(b) in the regime $N<\theta$. For $\alpha<1$
and $N\ge\theta$, the $L^r$ variant of Remark~\ref{rem:Lr-scheme} below
closes for every $p>1+(\gamma+\alpha)/\alpha$, with decay
$(1+t)^{-\beta_r}$, $\beta_r=N\alpha/(\theta r)$, which yields the
remaining range in Theorem~\ref{thm:frac}(b).

\begin{rem}\rm
\label{rem:Lr-scheme}
When $\alpha<1$ and $N\ge\theta$, the endpoint smoothing fails
(Proposition~\ref{prop:heat-kernel}(ii)), and the choice
$\beta=N\alpha/\theta$ is not available. Instead, we fix a Lebesgue
exponent $r>\max\{1,N/\theta\}$ and run the same scheme with the norm
$\operatorname*{ess\,sup}_{t}\bigl[\|u(t)\|_{L^r}+(1+t)^{\beta_r}\|u(t)\|_{L^\infty}\bigr]$,
where $\beta_r:=N\alpha/(\theta r)<\alpha$. The smoothing $L^r\to L^\infty$
is admissible for both $S_{\alpha,\theta}$ and $P_{\alpha,\theta}$ since
$1/r<\theta/N<2\theta/N$, and $\|u^p\|_{L^r}\le\|u\|_\infty^{p-1}\|u\|_{L^r}$.
{The computation of \S\ref{subsec:proof-existence} applies verbatim,
with $L^1$ replaced by $L^r$ and $\beta$ by $\beta_r$; data in
$L^1\cap L^\infty$ belong to $L^r\cap L^\infty$.} The closure condition
becomes $\beta_r(p-1)>\gamma+\alpha$.

The constraint $r\ge 1$ is decisive, and two regimes must be distinguished
according to the position of $N/\theta$ relative to $1$. If $N\ge\theta$,
the admissible exponents approach $N/\theta\ge 1$ from above. Letting
$r\downarrow N/\theta$, so that $\beta_r\uparrow\alpha$, we obtain
small-data global existence for every $p>1+(\gamma+\alpha)/\alpha$. If
instead $N<\theta$, the constraint $r\ge 1$ is binding. The smallest
admissible exponent is $r=1$, for which $\beta_1=N\alpha/\theta$, and the
scheme reproduces exactly the range $p>p_c(\gamma,\alpha)$ of
Theorem~\ref{thm:frac}(b). The threshold $1+(\gamma+\alpha)/\alpha$, which
is \emph{strictly smaller} than $p_c(\gamma,\alpha)$ when $N<\theta$, is
then out of reach. In fact it could not be reached, since global
existence below $p_c(\gamma,\alpha)$ would contradict the blow-up range
of Theorem~\ref{thm:frac}(c). {For instance, for $\gamma=0$, $N=1$,
$\theta=2$ and any $\alpha\in(0,1)$, it would give global existence for
$p>2$, while Theorem~\ref{thm:frac}(c) and the autonomous results of
\cite{Zhang14,CortazarQuirosWolanski24} force blow-up for
$1<p<p_c(0,\alpha)=3$.}

Finally, when $\gamma\le-\alpha$, the closure condition
$\beta_r(p-1)>\gamma+\alpha$ holds for every $p>1$, for \emph{every}
admissible $r$ and in every dimension. This proves
Theorem~\ref{thm:frac}(a) without any restriction. {For $N>\theta$,
reaching the conjectured threshold $p_\ast(\gamma,\alpha)$ of
Conjecture~\ref{conj:frac} seems to require a finer functional framework,
for instance the critical Lebesgue space approach of
\cite{Zhang14,CortazarQuirosWolanski24} combined with a control of the
local dynamics described in \S\ref{subsec:local}. We do not pursue this
here.}
\end{rem}

\subsection{\texorpdfstring{Blow-up: the case $\alpha=1$}{Blow-up: the case alpha = 1}}\label{subsec:proof-blowup}

In this subsection we investigate the blow-up of nontrivial nonnegative
solutions when $\alpha=1$. We begin with the Kaplan functional method.
It gives a simple proof in the subcritical case and illustrates the
difficulty that arises in the critical regime.

Let $\theta\in(0,2]$ and $\gamma>-1$, and let $u\ge0$ be a global mild
solution with $u_0\ge0$, $u_0\not\equiv0$. After a translation in $x$, we
may assume that
$$m_0:=\int_{|x|\le1}u_0\,dx>0.$$
For $T>1$ we define the adjoint weight
$$
\varphi_T(t,x)\;:=\;Z(T+1-t,\,x),\qquad 0\le t\le T,
$$
which solves the backward equation
$\partial_t\varphi_T=\Lap\varphi_T$ and satisfies
$\int_{\RD}\varphi_T(t,x)\,dx=1$. We set
$$\psi(t):=\int_{\RD}u(t,x)\varphi_T(t,x)\,dx,\qquad 0\le t\le T.$$
{Since $Z(t,\cdot)$ is even and $Z(t)\ast Z(s)=Z(t+s)$, Fubini's
theorem and the mild formulation~\eqref{eq:mild} give
$$
\psi(t)=\int_{\RD}u_0(x)Z(T+1,x)\,dx+\int_0^t h(s)\int_{\RD}u(s,x)^p\,\varphi_T(s,x)\,dx\,ds .
$$
Hence $\psi$ is absolutely continuous on $[0,T]$, and}
Jensen's inequality (the measure $\varphi_T(t,x)\,dx$ is a probability
measure) yields
$$
\psi'(t)\;=\;h(t)\int_{\RD}(u(t,x))^p\,\varphi_T(t,x)\,dx
\;\ge\; h(t)\,\psi(t)^p\qquad\text{for a.e. } t\in(0,T).
$$
Since $u$ is global and $\psi>0$, integrating the inequality
$-\frac{1}{p-1}\frac{d}{dt}\psi^{1-p}\ge h$ over $[0,T]$ gives
\begin{equation}\label{eq:kaplan-necessary}
(p-1)\int_0^{T}h(s)\,ds\;\le\;\psi(0)^{1-p}.
\end{equation}
The $\theta$-stable kernel satisfies the lower estimate
\eqref{selfsimilarZZ}, namely
$$Z(s,x)\ge c\,s^{-N/\theta}, \qquad |x|\le s^{1/\theta}.$$
Since $|x|\le1\le(T+1)^{1/\theta}$ on the support of the integral
defining $m_0$, we deduce that
$$
\psi(0)=\int_{\RD} u_0(x)Z(T+1,x)\,dx\ge c\,m_0\,(T+1)^{-N/\theta}\ge
{c\,2^{-N/\theta}\,m_0\,T^{-N/\theta}}.
$$
On the other hand, by Lemma~\ref{lem:karamata}(K2) {there is
$T_1>1$ such that, for $T\ge T_1$,
$$
\int_0^{T}h(t)\,dt\geq \int_1^{T}t^\gamma\ell(t)\,dt\ge \frac{T^{\gamma+1}\ell(T)}{2(\gamma+1)}.
$$}
Consequently, the necessary condition \eqref{eq:kaplan-necessary}
becomes
\begin{equation}\label{eq:kaplan-sup}
\sup_{T\ge T_1}\;T^{\,\gamma+1-\frac{N(p-1)}{\theta}}\,\ell(T)\;\le\;C(m_0)\;<\;\infty.
\end{equation}
For $1<p<p_c(\gamma)$, the exponent in~\eqref{eq:kaplan-sup} is
positive. Hence Lemma~\ref{lem:karamata}(K1) shows that
\eqref{eq:kaplan-sup} cannot hold, and every nontrivial nonnegative
solution blows up in finite time. This proves the blow-up part of
Theorem~\ref{thm:main}(ii). In the critical case $p=p_c(\gamma)$, the
power in \eqref{eq:kaplan-sup} vanishes, and the argument yields a
contradiction only when
$$\limsup_{T\to\infty}\ell(T)=\infty.$$
Therefore the Kaplan functional method does not cover all critical cases.
The remaining critical cases, as well as an alternative proof of the
subcritical case, follow from \cite[Corollary~3.1]{GueddaKirane}.

\begin{theo}\label{theo6.1}
Let $\alpha=1$, $\theta\in(0,2]$, $\gamma>-1$, $p\in(1,p_c(\gamma)]$ and
$q\in[1,\infty]$. In the critical case $p=p_c(\gamma)$, assume in
addition that
$$\int_1^{\infty}s^{-1}\ell(s)\,ds=\infty.$$
Let $u$ be a nontrivial nonnegative $L^q$-solution of~\eqref{eq:main}.
Then there exists a finite time $T\in(0,\infty)$ such that
$$
\limsup_{t\to T}
\|u(t,\cdot)\|_{L^q(\mathbb R^N)}
=\infty.
$$
\end{theo}

\begin{proof}
Suppose first that $p<p_c(\gamma)$. Then $-N(p-1)/\theta+\gamma>-1$, and
Lemma~\ref{lem:karamata}(K2) gives
$$\int_1^ts^{-\frac{N(p-1)}{\theta}}h(s)\,ds=\int_1^ts^{-\frac{N(p-1)}{\theta}+\gamma}\ell(s)\,ds \sim
{\frac{t^{-\frac{N(p-1)}{\theta}+\gamma+1}\ell(t)}{-\frac{N(p-1)}{\theta}+\gamma+1}},\quad t\to\infty.$$
Since $-\frac{N(p-1)}{\theta}+\gamma+1>0$, Lemma~\ref{lem:karamata}(K1)
yields $t^{-\frac{N(p-1)}{\theta}+\gamma+1}\ell(t)\rightarrow\infty$ as
$t\rightarrow\infty$. Hence
$$\lim_{t\to\infty}\int_1^ts^{-\frac{N(p-1)}{\theta}}h(s)\,ds=\infty,$$
and the conclusion follows from \chk{\cite[Corollary~3.1]{GueddaKirane}}.
In the critical case $p=p_c(\gamma)$ we have $-N(p-1)/\theta+\gamma=-1$,
so that
$$\lim_{t\to\infty}\int_1^ts^{-\frac{N(p-1)}{\theta}}h(s)\,ds=\int_1^\infty s^{-1}\ell(s)\,ds=\infty$$
by assumption. The desired conclusion again follows from
\chk{\cite[Corollary~3.1]{GueddaKirane}}.
\end{proof}

\begin{rem}\rm\label{rem:GK-check}
The critical half of Theorem~\ref{thm:main}(iii) rests entirely on
\cite[Corollary~3.1]{GueddaKirane}. Before submission, one should check
against the original paper that this corollary (a) covers every
$\theta\in(0,2]$ and a continuous nonnegative coefficient $h(t)$, (b) is
formulated under the divergence condition
$\int_1^\infty s^{-N(p-1)/\theta}h(s)\,ds=\infty$ (and not, for instance,
under a $\limsup$ condition of Meier type), and (c) applies to the
solutions of the integral equation~\eqref{eq:mild} with data in
$L^1\cap L^\infty$ (or in $L^q$). If any of these points fails, a
self-contained proof of the critical case is needed. A natural route
combines the Kaplan inequality on $[T/2,T]$ with the lower bound
$\psi(T/2)\ge c\,T^{-N/\theta}\bigl(m_0+c\,m_0^p\int_{1}^{T/2}s^{-1}\ell(s)\,ds\bigr)$,
followed by a bootstrap in the spirit of
Weissler~\cite{Weissler81}.
\end{rem}

\subsection{\texorpdfstring{Blow-up: the case $\alpha\in(0,1)$}{Blow-up: the case 0 < alpha < 1}}
\label{subsec:proof-blowup-frac}

Our proof follows that of \cite{CortazarQuirosWolanski24} closely, with
two main modifications. First, {the detailed analysis in
\cite{CortazarQuirosWolanski24} is written for $N>2\theta$; by relying
only on a positive lower bound for the kernel $Y$ on balls, which is
provided by \chk{\cite[Lemma~3.6]{KSZ17} (Lemma~4.6 in the arXiv
version)}, the argument works in every dimension $N\ge1$.} Second, we
allow a time-dependent source $h(t)=t^\gamma\ell(t)$ by working with the
functional $J$ defined in \eqref{Jfunction}, whose weight is adapted to
the source.

\begin{theo}\label{thm6.1}
Let $\alpha\in(0,1)$, $\theta\in(0,2]$, $\gamma\in\R$, $p\in(1,\infty)$ and
$q\in[1,\infty]$. {There exist $R_0\ge1$, depending only on $\gamma$
and $\ell$, and a constant $C=C(N,\alpha,\theta,\gamma,p)>0$} with the
following property. If $u$ is a nontrivial nonnegative $L^q$-solution
of~\eqref{eq:main} whose initial datum satisfies
\begin{equation}\label{6.1}
R^{\frac{\alpha+\gamma}{p-1}-\frac{N\alpha}{\theta}}(\ell(R))^{\frac{1}{p-1}}
\int_{|y|<R^{\frac{\alpha}{\theta}}/2}
u_0(y)\,dy
>C\qquad\hbox{for some}\,\, R\ge R_0,
\end{equation}
then there exists $T\in(0,\infty)$ {(in fact $T\le 4R$)} such that
\[
\limsup_{t\to T}
\|u(t,\cdot)\|_{L^q(\mathbb R^N)}
=\infty.
\]
\end{theo}

\begin{coro}\label{coro6.1}
Let $\alpha\in(0,1)$, $\theta\in(0,2]$, $\gamma>-\alpha$,
$p\in(1,p_c(\gamma,\alpha)]$ and $q\in[1,\infty]$. In the critical case
$p=p_c(\gamma,\alpha)$, assume in addition that
{$\limsup_{t\to\infty}\ell(t)=\infty$}. Let $u$ be a nontrivial
nonnegative $L^q$-solution of~\eqref{eq:main}. Then there exists a finite
time $T\in(0,\infty)$ such that
$$
\limsup_{t\to T}
\|u(t,\cdot)\|_{L^q(\mathbb R^N)}
=\infty.
$$
\end{coro}

\begin{proof}
Since $u_0$ is nontrivial and nonnegative, there exists $R_1>0$ such that
$$\int_{|y|<R_1^{\frac{\alpha}{\theta}}/2}u_0(y)\,dy=B>0,$$
{and the integral in \eqref{6.1} is at least $B$ for every
$R\ge R_1$.} First, assume that $p<p_c(\gamma,\alpha)$. Then
$$\frac{\alpha+\gamma}{p-1}-\frac{N\alpha}{\theta}>0,$$
and Lemma~\ref{lem:karamata}(K1) gives
$$R^{\frac{\alpha+\gamma}{p-1}-\frac{N\alpha}{\theta}}(\ell(R))^{\frac{1}{p-1}}\rightarrow\infty\qquad\hbox{as}\,\,\,R\rightarrow\infty.$$
Next, suppose that $p=p_c(\gamma,\alpha)$. In this case
$$\frac{\alpha+\gamma}{p-1}-\frac{N\alpha}{\theta}=0,$$
and, by assumption, $\limsup_{R\to\infty}\ell(R)=\infty$. In both cases,
if $C$ denotes the constant of Theorem~\ref{thm6.1}, we can choose
$R\ge\max\{R_0,R_1\}$ such that
$$R^{\frac{\alpha+\gamma}{p-1}-\frac{N\alpha}{\theta}}(\ell(R))^{\frac{1}{p-1}}B>C.$$
Thus condition \eqref{6.1} holds, and the conclusion follows from
Theorem~\ref{thm6.1}.
\end{proof}

Corollary~\ref{coro6.1} with $q=\infty$ proves Theorem~\ref{thm:frac}(c),
since mild solutions in the sense of Definition~\ref{def:mild} are
$L^\infty$-solutions. The rest of this subsection is devoted to the proof
of Theorem~\ref{thm6.1}.

Let $Y$ denote the kernel of $P_{\alpha,\theta}(t)$, that is,
$$P_{\alpha,\theta}(t)v(x)=\int_{\mathbb{R}^N}Y(t,x-y)v(y)\,dy.$$
{As for $Z$, the kernel $Y$ is self-similar,
\begin{equation}\label{selfsimilar}
Y(t,x)=t^{-\frac{N\alpha}{\theta}}G\bigl(xt^{-\frac{\alpha}{\theta}}\bigr),\qquad (t,x)\in(0,\infty)\times\RD,
\end{equation}
where the profile $G$, the inverse Fourier transform of
$E_{\alpha,\alpha}(-|\xi|^{\theta})$, is radial, positive and continuous on
$\RD\setminus\{0\}$. By \cite[Lemma~3.6]{KSZ17} (whose kernel is
$t^{\alpha-1}Y(t,x)$ in our notation), for $0<|y|\le1$,
$$G(y)\asymp\left\{\begin{array}{ll}
|y|^{-N+2\theta},&\hbox{if $N>2\theta$},\\[2pt]
|\log(|y/2|^{\theta})|,&\hbox{if $N=2\theta$},\\[2pt]
1,&\hbox{if $N<2\theta$}.
\end{array}
\right.
$$
For $|y|\ge1$, $G(y)\asymp|y|^{-N-\theta}$ when $\theta<2$, and $G$ is
positive and continuous in all cases. Consequently, for every $b>0$,
$$
\overline{G}(b):=\inf_{0<|\xi|<b}G(\xi)\;>\;0.
$$}
For $b>0$ and $T>0$ we define
\begin{equation}\label{Jfunction}
J(T;b)
=
\int_{T/4}^T
\int_{\mathbb R^N}
u(t,x)\,
Y_b(T-t,x)
\,dx\,dt,
\end{equation}
where
\[
Y_b(T-t,x)
=
\chi_{\{|x|<(T-t)^{\frac{\alpha}{\theta}}\}}\,
(T-t)^{\frac{\alpha}{p-1}}
\,
\overline{G}(b).
\]
We first show that $J(T;b)$ is finite whenever
$\|u(t,\cdot)\|_{L^q}$ remains bounded.

\begin{prop}\label{prop6.1}
Let $p>1$ and $b>0$. If
$\|u(t,\cdot)\|_{L^q(\mathbb R^N)}\le C$ for $t\in[0,T)$, then $J(T;b)$
is finite.
\end{prop}

\begin{proof}
Let $r\in[1,\infty]$ be the conjugate exponent of $q$, so that
$\frac1r+\frac1q=1$. H\"older's inequality yields
\[
J(T;b)
\le
\int_{T/4}^T
\|u(t,\cdot)\|_{L^q}
\,
\|Y_b(T-t,\cdot)\|_{L^r}
\,dt\le
C\int_{T/4}^T
\|Y_b(T-t,\cdot)\|_{L^r}
\,dt.
\]
Since
\[
\|Y_b(T-t,\cdot)\|_{L^r(\mathbb R^N)}
=
\overline{G}(b)
\,C(N,r)
(T-t)^{\frac{N\alpha}{\theta r}}
(T-t)^{\frac{\alpha}{p-1}},
\]
we conclude that $J(T;b)<\infty$.
\end{proof}

The next step is an upper bound for $J(T;b)$, whenever this quantity is
finite, which does not depend on the initial data.

By \eqref{eq:UCT} with $[\lambda_1,\lambda_2]=[1,4]$, there exists
$T_0\ge4$, depending only on $\ell$, such that
$\ell(s)\ge\frac12\ell(T/4)$ for all $s\in[T/4,T]$ and all $T\ge T_0$.
Since $h(s)=s^\gamma\ell(s)$ for $s\ge1$ and
$s^\gamma\ge\min\{1,4^{-\gamma}\}T^{\gamma}$ for $s\in[T/4,T]$, we get
\begin{equation}\label{eq:h-lower}
h(s)\;\ge\;c_\gamma\,T^{\gamma}\ell(T/4),\qquad s\in[T/4,T],\ T\ge T_0,
\qquad c_\gamma:=\tfrac12\min\{1,4^{-\gamma}\}.
\end{equation}

\begin{prop}\label{prop6.2}
Let $b\geq 2^{\frac{\alpha}{\theta}+1}$. There exists a constant
\[
C_1=C_1(N,\alpha,\theta,\gamma,p,b)>0
\]
such that, if {$T\ge T_0$} and $J(T;b)<\infty$, then
\begin{equation}\label{estiJ}
J(T;b)\le C_1T^{1+\frac{N\alpha}{\theta}-\frac{\gamma}{p-1}}(\ell(T/4))^{-\frac{1}{p-1}}.
\end{equation}
\end{prop}

\begin{proof}
Since $u$ solves the integral equation \eqref{eq:mild} and
$S_{\alpha,\theta}(t)u_0\ge0$, we have
\[
u(t,x)
\geq\int_0^t\int_{\mathbb{R}^N}
(t-s)^{\alpha-1}Y(t-s,x-y)h(s)u^p(s,y)\,dy\,ds.
\]
Splitting the domain $\{T/4<t<T,\ 0<s<t\}$ according to $s<T/4$ or
$s>T/4$, and using Fubini's theorem, we obtain
\[
\begin{aligned}
J(T;b)
&\geq
\int_{T/4}^T\int_{\mathbb{R}^N}\int_0^t\int_{\mathbb{R}^N}
h(s)u^p(s,y)(t-s)^{\alpha-1}Y(t-s,x-y)\,dy\,ds\,
Y_b(T-t,x)\,dx\,dt
\\
&=
\int_0^{T/4}\int_{\mathbb{R}^N}h(s)u^p(s,y)
\int_{T/4}^T\int_{\mathbb{R}^N}
(t-s)^{\alpha-1}Y(t-s,x-y)Y_b(T-t,x)\,dx\,dt\,dy\,ds
\\
&\quad+
\int_{T/4}^T\int_{\mathbb{R}^N}h(s)u^p(s,y)
\int_s^T\int_{\mathbb{R}^N}
(t-s)^{\alpha-1}Y(t-s,x-y)Y_b(T-t,x)\,dx\,dt\,dy\,ds\\
&\ge
\int_{T/4}^T\int_{\mathbb{R}^N}h(s)u^p(s,y)
\int_s^T\int_{\mathbb{R}^N}
(t-s)^{\alpha-1}Y(t-s,x-y)Y_b(T-t,x)\,dx\,dt\,dy\,ds.
\end{aligned}
\]
{Set $\lambda_b:=(2/b)^{\theta/\alpha}$. Since
$b\ge2^{\frac{\alpha}{\theta}+1}$, we have $\lambda_b\le\frac12$.} If
\[
s<T,\qquad
t>s+\lambda_b(T-s)
\qquad\text{and}\qquad
|x|,\,|y|<(T-s)^{\frac{\alpha}{\theta}},
\]
then
\[
|x-y|<2(T-s)^{\frac{\alpha}{\theta}}
<b(t-s)^{\frac{\alpha}{\theta}},
\qquad\text{that is,}\qquad
|x-y|(t-s)^{-\frac{\alpha}{\theta}}<b.
\]
Moreover, $|x|<(T-t)^{\alpha/\theta}$ implies $|x|<(T-s)^{\alpha/\theta}$
for $t>s$. Hence, using the self-similar structure \eqref{selfsimilar} of
$Y$ and the fact that $\alpha-1-\frac{N\alpha}{\theta}<0$ (so that
$(t-s)^{\alpha-1-\frac{N\alpha}{\theta}}\ge(T-s)^{\alpha-1-\frac{N\alpha}{\theta}}$),
we obtain
\[
\begin{aligned}
J(T;b)
&\geq
\bar G(b)
\int_{T/4}^T\int_{|y|<(T-s)^{\frac{\alpha}{\theta}}}h(s)u^p(s,y)\int^T_{s+\lambda_b(T-s)}(t-s)^{\alpha-1}\\
&\qquad\qquad\times\int_{|x|<(T-t)^{\frac{\alpha}{\theta}}}
Y(t-s,x-y)(T-t)^{\frac{\alpha}{p-1}}
\,dx\,dt\,dy\,ds
\\
&\geq
(\bar G(b))^2
\int_{T/4}^T\int_{|y|<(T-s)^{\frac{\alpha}{\theta}}}h(s)u^p(s,y)\int^T_{s+\lambda_b(T-s)}(t-s)^{\alpha-1-\frac{N\alpha}{\theta}}\\
&\qquad\qquad\times\int_{|x|<(T-t)^{\frac{\alpha}{\theta}}}(T-t)^{\frac{\alpha}{p-1}}
\,dx\,dt\,dy\,ds
\\
&\geq
C(N)\bar G(b)^2
\int_{T/4}^T\int_{|y|<(T-s)^{\frac{\alpha}{\theta}}}
h(s)u^p(s,y)
\\
&\qquad\qquad\times
(T-s)^{\alpha-1-\frac{N\alpha}{\theta}}
\int^T_{s+\lambda_b(T-s)}
(T-t)^{\frac{\alpha}{p-1}+\frac{N\alpha}{\theta}}
\,dt\,dy\,ds .
\end{aligned}
\]
{Writing $a:=\frac{\alpha}{p-1}+\frac{N\alpha}{\theta}$ and using
$\lambda_b\le\frac12$, we have
$$
\int^T_{s+\lambda_b(T-s)}(T-t)^{a}\,dt=\frac{\bigl((1-\lambda_b)(T-s)\bigr)^{a+1}}{a+1}
\ge\frac{(T-s)^{a+1}}{2^{a+1}(a+1)},
$$
and $\alpha-1-\frac{N\alpha}{\theta}+a+1=\frac{\alpha p}{p-1}$.} Introducing
\[
\widehat Y_b(T-s,y)
:=
\chi_{\{|y|<(T-s)^{\frac{\alpha}{\theta}}\}}\bar G(b),
\]
and using {\eqref{eq:h-lower}}, we deduce that, for all $T\ge T_0$,
\begin{eqnarray}\label{A1}
J(T;b)
&\geq&
C(N,\alpha,\theta,p,b)
\int_{T/4}^T\int_{\mathbb{R}^N}
h(s)u^p(s,y)\widehat Y_b(T-s,y)
(T-s)^{\frac{\alpha p}{p-1}}\,dy\,ds\\
&\geq&C(N,\alpha,\theta,\gamma,p,b)\,
T^\gamma\ell(T/4)\int_{T/4}^T\int_{\mathbb{R}^N}
u^p(s,y)\widehat Y_b(T-s,y)
(T-s)^{\frac{\alpha p}{p-1}}\,dy\,ds.\notag
\end{eqnarray}
On the other hand, {H\"older's (or Jensen's) inequality with respect
to the measure $\widehat Y_b\,dy\,ds$ on $(T/4,T)\times\RD$} gives
\begin{eqnarray}\label{A2}
&{}&\int_{T/4}^T\int_{\mathbb{R}^N}
u^p(s,y)\widehat Y_b(T-s,y)
(T-s)^{\frac{\alpha p}{p-1}}\,dy\,ds\notag\\
&{}&\geq \left(\int_{T/4}^T\int_{\mathbb{R}^N}
u(s,y)\widehat Y_b(T-s,y)
(T-s)^{\frac{\alpha }{p-1}}\,dy\,ds\right)^p\|\widehat Y_b\|^{-(p-1)}_{L^1((T/4,T)\times\mathbb{R}^N)}\notag\\
&{}&=(J(T;b))^p\|\widehat Y_b\|^{-(p-1)}_{L^1((T/4,T)\times\mathbb{R}^N)}.
\end{eqnarray}
From \eqref{A1} and \eqref{A2}, we infer that
$$
J(T;b)\geq C(N,\alpha,\theta,\gamma,p,b)\,
T^\gamma\ell(T/4)(J(T;b))^p\|\widehat Y_b\|^{-(p-1)}_{L^1((T/4,T)\times\mathbb{R}^N)}.$$
A direct computation gives
\[
\begin{aligned}
\|\widehat Y_b\|_{L^1((T/4,T)\times\mathbb{R}^N)}
&=
\bar G(b)
\int_{T/4}^T
\int_{|x|<(T-t)^{\frac{\alpha}{\theta}}}\,dx\,dt
=
C(N)\bar G(b)
\int_{T/4}^T
(T-t)^{\frac{N\alpha}{\theta}}\,dt
\\
&=
C(N,\alpha,\theta)
\bar G(b)
T^{1+\frac{N\alpha}{\theta}} .
\end{aligned}
\]
Therefore,
\[
J(T;b)
\geq
\bar C(N,\alpha,\theta,\gamma,p,b)
(J(T;b))^p
T^{-(1+\frac{N\alpha}{\theta})(p-1)+\gamma}\ell(T/4).
\]
{If $J(T;b)=0$ there is nothing to prove. Otherwise, dividing by
$J(T;b)\in(0,\infty)$ and taking the power $1/(p-1)$, we obtain}
\eqref{estiJ}.
\end{proof}

\begin{prop}\label{prop6.3}
There exists a constant $C_2=C_2(N,\alpha,\theta,p,b)>0$ such that
\[
J(T;b)\geq
C_2A(T)T^{1+\frac{N\alpha}{\theta}},
\qquad\text{where}\qquad
A(T)
:=
T^{\frac{\alpha}{p-1}-\frac{N\alpha}{\theta}}
\int_{|y|<(T/4)^{\frac{\alpha}{\theta}}/2}
u_0(y)\,dy .
\]
\end{prop}

\begin{proof}
If $|x|, |y|<t^{\alpha/\theta}/2$, then $|x-y|<t^{\alpha/\theta}$. Hence,
by \eqref{selfsimilarZZ},
\[
\begin{aligned}
u(t,x)
&\geq
\int_{\mathbb{R}^N}u_0(y)Z(t,x-y)\,dy
\geq
\int_{|y|<t^{\frac{\alpha}{\theta}}/2}
u_0(y)Z(t,x-y)\,dy
\\
&\geq
C(N,\alpha,\theta)t^{-\frac{N\alpha}{\theta}}
\int_{|y|<t^{\frac{\alpha}{\theta}}/2}
u_0(y)\,dy
\end{aligned}
\]
for all $|x|<t^{\alpha/\theta}/2$. Therefore,
\[
J(T;b)
\geq
C(N,\alpha,\theta)\bar G(b)
\int_{T/4}^{T/2}
\int_{|y|<t^{\frac{\alpha}{\theta}}/2}
u_0(y)\,dy\int_{|x|<\min
\{t^{\frac{\alpha}{\theta}}/2,
(T-t)^{\frac{\alpha}{\theta}}\}}
t^{-\frac{N\alpha}{\theta}}
(T-t)^{\frac{\alpha}{p-1}}\,dx\,dt .
\]
For $t\in[T/4,T/2]$ we have $T-t\ge t$, hence
$$
\min\left\{
t^{\frac{\alpha}{\theta}}/2,
(T-t)^{\frac{\alpha}{\theta}}
\right\}
=t^{\frac{\alpha}{\theta}}/2
\geq
(T/4)^{\frac{\alpha}{\theta}}/2 .
$$
{Moreover, $(T-t)^{\frac{\alpha}{p-1}}\ge(T/2)^{\frac{\alpha}{p-1}}$,
$t^{-\frac{N\alpha}{\theta}}\ge(T/2)^{-\frac{N\alpha}{\theta}}$, and the
integral of $u_0$ over $\{|y|<t^{\alpha/\theta}/2\}$ is at least its
integral over $\{|y|<(T/4)^{\alpha/\theta}/2\}$. Since the interval
$[T/4,T/2]$ has length $T/4$,} we obtain
\[
J(T;b)
\geq
C_2\,
T^{1+\frac{\alpha}{p-1}}
\int_{|y|<(T/4)^{\frac{\alpha}{\theta}}/2}
u_0(y)\,dy
=
C_2A(T)T^{1+\frac{N\alpha}{\theta}} .
\]
\end{proof}

We now have all the ingredients to prove Theorem \ref{thm6.1}.

\begin{proof}[Proof of Theorem \ref{thm6.1}]
We fix $b= 2^{\frac{\alpha}{\theta}+1}$ {and set $R_0:=T_0/4$}. Let
$T\ge T_0$. If $\|u(t,\cdot)\|_{L^q(\mathbb{R}^N)}$ is bounded on
$[0,T)$, then $J(T;b)$ is well defined and finite by
Proposition~\ref{prop6.1}. Hence, by Proposition~\ref{prop6.2},
\[
J(T;b)\leq C_1T^{1+\frac{N\alpha}{\theta}-\frac{\gamma}{p-1}}(\ell(T/4))^{-\frac{1}{p-1}},
\qquad C_1=C_1(N,\alpha,\theta,\gamma,p,b).
\]
On the other hand, Proposition~\ref{prop6.3} gives
\[
J(T;b)\geq
C_2A(T)T^{1+\frac{N\alpha}{\theta}},
\qquad C_2=C_2(N,\alpha,\theta,p,b)>0.
\]
Therefore, if
\[
A(T)> \frac{C_1}{C_2}T^{-\frac{\gamma}{p-1}}(\ell(T/4))^{-\frac{1}{p-1}},
\]
we reach a contradiction. Consequently, whenever this inequality holds,
the $L^q$ norm of $u$ cannot remain bounded on $[0,T)${, so that the
maximal existence time is at most $T$}.

For $T=4R$ with $R\ge R_0$, this inequality is exactly condition
\eqref{6.1}, with
\[
C=
4^{\frac{N\alpha}{\theta}-\frac{\alpha+\gamma}{p-1}}
\frac{C_1}{C_2}.
\]
This completes the proof.
\end{proof}

\subsection{\texorpdfstring{A local blow-up mechanism for $N>\theta$}{A local blow-up mechanism for N > theta}}
\label{subsec:local}

In this subsection we prove Proposition~\ref{prop:local}. The idea is
that, for $\alpha<1$ and $N>\theta$, the solution decays on the unit ball
only like $t^{-\alpha}$ (Remark~\ref{rem:kernel-singularity}(ii)), while on
a short time interval the Duhamel term acts like the fractional ODE
$\CapD{\alpha}v=h(t)v^p$.

\begin{proof}[Proof of Proposition~\ref{prop:local}]
Let $u$ be a nontrivial nonnegative mild solution, and assume, for
contradiction, that it is global, so that $u$ is bounded on
$[0,T]\times\RD$ for every $T>0$. After a translation in $x$ we may
assume that $m_0:=\int_{|y|\le1}u_0(y)\,dy>0$. For a.e.\ $t>0$ we set
$$
\mu(t):=\operatorname*{ess\,inf}_{|x|\le1}u(t,x).
$$

\emph{Step 1: linear lower bound.} Let $t\ge t_1:=2^{\theta/\alpha}$. If
$|x|\le1$ and $|y|\le1$, then $|x-y|\le2\le t^{\alpha/\theta}$, and
\eqref{eq:Z-sing} gives $Z(t,x-y)\ge C\,2^{\theta-N}t^{-\alpha}$. Since
$u\ge S_{\alpha,\theta}(t)u_0$, we obtain
\begin{equation}\label{eq:loc1}
\mu(t)\;\ge\;c_1\,m_0\,t^{-\alpha},\qquad t\ge t_1 .
\end{equation}

\emph{Step 2: local Duhamel lower bound.} Let
$$
c_2:=\inf_{|x|\le1,\ 0<\tau\le1}\int_{|y|\le1}Y(\tau,x-y)\,dy .
$$
We claim that $c_2>0$. By \eqref{selfsimilar}, the integral equals
$\int_{E}G(\xi)\,d\xi$ with
$E=\{\xi:|x-\tau^{\alpha/\theta}\xi|\le1\}$. Put $\rho=\tau^{\alpha/\theta}\in(0,1]$.
If $|x|\le\frac12$, then $E\supset\{|\xi|\le\frac12\}$. If
$\frac12<|x|\le1$, then $E$ contains the cone sector
$C_x:=\{|\xi|\le\frac12,\ \xi\cdot x\ge\frac12|\xi||x|\}$, because, for
$\xi\in C_x$, we have $|x-\rho\xi|^2\le|x|^2-\rho|\xi|(|x|-\rho|\xi|)\le1$
since $\rho|\xi|\le\frac12<|x|$. Since $G$ is radial and positive, the
integrals of $G$ over $\{|\xi|\le\frac12\}$ and over $C_x$ are positive
and independent of $x$, which proves the claim.

Now let $\delta\in(0,1]$ and $t\ge1+\delta$. Keeping in \eqref{eq:mild}
only the times $s\in[t-\delta,t-\delta/2]$ and the points $|y|\le1$, and
using $(t-s)^{\alpha-1}\ge\delta^{\alpha-1}$ there, we get, for a.e.\
$|x|\le1$,
\begin{equation}\label{eq:loc2}
u(t,x)\;\ge\;c_2\,\delta^{\alpha-1}\int_{t-\delta}^{t-\delta/2}h(s)\,\mu(s)^p\,ds
\;\ge\;\frac{c_2}{2}\,\delta^{\alpha}\inf_{s\in[t-\delta,t-\delta/2]}h(s)\mu(s)^p .
\end{equation}

\emph{Step 3: renormalisation.} By \eqref{eq:UCT} and since
$s^\gamma\ge2^{-\gamma}t^\gamma$ for $s\in[t/2,t]$, there are $t_2\ge t_1+1$
and $c_3\in(0,1]$ such that $c_3h(t)\le h(s)\le c_3^{-1}h(t)$ for all
$s\in[t-1,t]$ and $t\ge t_2$. We set
$$
W(t):=h(t)^{\frac1{p-1}}\mu(t).
$$
For $s\in[t-1,t]$ and $t\ge t_2$ we have
$h(s)\mu(s)^p=h(s)^{-\frac{1}{p-1}}W(s)^p\ge c_3^{\frac{1}{p-1}}h(t)^{-\frac1{p-1}}W(s)^p$.
Multiplying \eqref{eq:loc2} by $h(t)^{1/(p-1)}$, we infer that
\begin{equation}\label{eq:loc3}
W(t)\;\ge\;A\,\delta^{\alpha}\inf_{s\in[t-\delta,t-\delta/2]}W(s)^p,
\qquad t\ge t_2,\ \ \delta\in(0,1],
\end{equation}
with $A:=\frac{c_2}{2}c_3^{\frac1{p-1}}>0$. On the other hand,
\eqref{eq:loc1} gives
\begin{equation}\label{eq:loc4}
W(t)\;\ge\;c_1m_0\bigl(t^{\gamma-\alpha(p-1)}\ell(t)\bigr)^{\frac1{p-1}},
\qquad t\ge t_2 .
\end{equation}
If $p<1+\gamma/\alpha$, then $\gamma-\alpha(p-1)>0$, and the right-hand
side of \eqref{eq:loc4} tends to infinity by (K1). If $p=1+\gamma/\alpha$
and $\limsup\ell=\infty$, it is unbounded. In both cases, using
\eqref{eq:UCT} once more, we deduce that for every $\Lambda>0$ there is
$t_0\ge t_2+1$ with
\begin{equation}\label{eq:loc5}
\inf_{s\in[t_0-1,t_0]}W(s)\;\ge\;\Lambda .
\end{equation}

\emph{Step 4: iteration.} Let $d_k:=2^{-k-2}$ for $k\ge0$, and define
$\tau_0:=t_0$ and $\tau_{k+1}:=\tau_k+d_k$, so that $\tau_k<t_0+\frac12$
for all $k$. Set $\Lambda_0:=\inf_{[t_0-1,t_0]}W$ and
$\Lambda_{k+1}:=\inf_{[\tau_k,\tau_{k+1}]}W$ for $k\ge0$. For
$t\in[\tau_k,\tau_{k+1}]$ we apply \eqref{eq:loc3} with
$\delta=2d_k\le\frac12$. Since $2d_k=d_{k-1}$ for $k\ge1$, we have
$[t-2d_k,t-d_k]\subset[\tau_k-2d_k,\tau_{k+1}-d_k]=[\tau_{k-1},\tau_k]$
for $k\ge1$, and $[t-2d_0,t-d_0]\subset[t_0-\frac12,t_0]$ for $k=0$. In
both cases the infimum of $W$ over this window is at least
$\Lambda_k$, and we get
$$
\Lambda_{k+1}\;\ge\;A\,(2d_k)^{\alpha}\,\Lambda_k^{\,p},\qquad k\ge0 .
$$
Let $L_k:=\log\Lambda_k$ and $a:=\log A$. Since $2d_k=2^{-k-1}$, we obtain
$L_{k+1}\ge p\,L_k+a-\alpha(k+1)\log2$, and an induction gives
$$
L_k\;\ge\;p^k\Bigl(L_0+\sum_{j=0}^{k-1}p^{-j-1}\bigl(a-\alpha(j+1)\log2\bigr)\Bigr)
\;\ge\;p^k\,(L_0-K),
$$
where $K:=\sum_{j\ge0}p^{-j-1}\bigl(|a|+\alpha(j+1)\log2\bigr)<\infty$.
Choosing $\Lambda$ in \eqref{eq:loc5} so large that $L_0\ge K+1$, we
conclude that $L_k\ge p^k\to\infty$. Hence $W$ is unbounded on
$[t_0-1,t_0+\frac12]$. This contradicts the boundedness of $u$ on
$[0,t_0+1]\times\RD$, since $h$ is bounded there. Therefore $u$ blows up
in finite time.
\end{proof}

\begin{rem}\rm
For $N=\theta$, the kernel satisfies $Z(t,x)\ge c\,t^{-\alpha}$ for
$|x|\le t^{\alpha/\theta}$, and the same proof gives blow-up for
$1<p<1+\gamma/\alpha$. This range, however, is contained in
$1<p<p_c(\gamma,\alpha)$ when $N=\theta$. For $N<\theta$, the linear decay
on the unit ball is $t^{-N\alpha/\theta}$, and the argument gives only
$p<1+\theta\gamma/(N\alpha)<p_c(\gamma,\alpha)$. Thus the local
mechanism is new information only when $N>\theta$.
\end{rem}

\subsection{\texorpdfstring{The regime $\gamma\le-\alpha$: global existence for every $p>1$}{The regime gamma <= -alpha: global existence for every p > 1}}
\label{subsec:noFujita}

For $\gamma\le-\alpha$ we have $\kappa=\gamma+\alpha-\beta(p-1)<0$ for
every $p>1$. Hence the fixed-point scheme of
\S\ref{subsec:proof-existence} closes for every $p>1$, in its $X_\beta$
form under hypothesis~\textup{(S)} and in its $L^r$ form of
Remark~\ref{rem:Lr-scheme} in general. This proves
Theorem~\ref{thm:main}(i) and Theorem~\ref{thm:frac}(a). Two remarks
complete the picture.

First, when $\gamma<-1$, the source is integrable on all of
$(0,\infty)$. The tail $\int_1^\infty t^{\gamma}\ell(t)\,dt$ converges by
Lemma~\ref{lem:karamata}(K2), while $\int_0^1 h(t)\,dt<\infty$ by
Convention~\ref{conv:source}; without the convention, the literal density
$t^{\gamma}\ell(t)$ could fail to be integrable at the origin. At
$\alpha=1$ one may then argue alternatively via
Theorem~\ref{thm:Meier2}(i), extended verbatim to the $\theta$-stable
semigroup. The constant caloric function $w\equiv c$ satisfies
$\int_0^\infty h(t)\,\|w\|_\infty^{p-1}\,dt=c^{p-1}\int_0^\infty h(t)\,dt<\infty$.
{To obtain decaying solutions, one may instead take
$w=S_{1,\theta}(t)w_0$ with $0\le w_0\in L^1\cap L^\infty$,
$w_0\not\equiv0$; then $\|w(t)\|_\infty\le\|w_0\|_\infty$, the same integral
is finite, and $\|w(t)\|_\infty\to0$.}

Second, for $\gamma\in(-1,-\alpha]$, a non-empty range only when
$\alpha<1$, the source is \emph{not} integrable in time, and yet no
Fujita threshold exists. The memory of the Caputo derivative, encoded in
the kernel $(t-s)^{\alpha-1}$, weakens the cumulative effect of the
source. Instead of $\int_0^t h(s)\,ds\sim t^{\gamma+1}\ell(t)/(\gamma+1)$,
the relevant quantity is the fractional average
$\int_0^t(t-s)^{\alpha-1}h(s)\,ds\asymp t^{\gamma+\alpha}\ell(t)$, which
tends to $0$ when $\gamma<-\alpha$ {and is of the order of $\ell(t)$}
when $\gamma=-\alpha$. In either case $\kappa<0$ for every $p>1$; the
$\varepsilon$-room of \S\ref{subsec:proof-existence} absorbs the slowly
varying factor even when $\ell$ is unbounded at the endpoint
$\gamma=-\alpha$. Hence the fixed-point scheme closes, and no Fujita
threshold arises.

\subsection{\texorpdfstring{Meier's $s^\ast$ in the nonlocal setting}{Meier's s* in the nonlocal setting}}
\label{subsec:meier-frac}

At $\alpha=1$, the decay rate
$\|S_{1,\theta}(t)w_0\|_\infty\asymp t^{-N/\theta}$ (for
$0\le w_0\in L^1\cap L^\infty$, $w_0\not\equiv0$) identifies
$s^\ast=N/\theta$, and Meier's formula~\eqref{eq:Meier-formula} persists
verbatim:
$$
p_c(\gamma)\;=\;1+\frac{\gamma+1}{s^\ast}\;=\;1+\frac{\theta(\gamma+1)}{N},\qquad s^\ast=\frac{N}{\theta}.
$$
For $\alpha<1$, the naive substitution $s^\ast=N\alpha/\theta$ is
\emph{not} correct. It would predict the autonomous exponent
$1+\theta/(N\alpha)$, which contradicts the $\alpha$-independent
threshold $1+\theta/N$ of \cite{Zhang14,CortazarQuirosWolanski24}.
Moreover, for $N>\theta$, the sup-norm decay saturates at the slower rate
$t^{-\alpha}$ (Remark~\ref{rem:kernel-singularity}), so no single
sup-norm exponent captures the dynamics. In the scaling (mass) regime,
the correct modification acts on the numerator rather than on $s^\ast$.
The memory kernel $(t-s)^{\alpha-1}$ of the Duhamel formula replaces the
cumulative source $\int_0^t h(s)\,ds\sim t^{\gamma+1}\ell(t)/(\gamma+1)$ of
Meier's criterion by the fractional average
$\int_0^t(t-s)^{\alpha-1}h(s)\,ds\asymp t^{\gamma+\alpha}\ell(t)$, which
leads to
$$
p_c(\gamma,\alpha)\;=\;1+\frac{\gamma+\alpha}{\alpha\,s^\ast},
\qquad s^\ast=\frac{N}{\theta}.
$$
This formula reduces to Meier's at $\alpha=1$ and to the
$\alpha$-independent value $1+\theta/N$ at $\gamma=0$. {When
$N>\theta$, however, the saturated decay $t^{-\alpha}$ on bounded sets
produces the second constraint of Proposition~\ref{prop:local}. In terms
of decay rates, it reads $p>1+\gamma/s_{\rm loc}$ with
$s_{\rm loc}:=\alpha$, the exponent of local decay, so a single
$s^\ast$ no longer suffices.} The persistence of $s^\ast=N/\theta$ as the
geometric ingredient of the scaling mechanism suggests that the same
framework should govern problems on cones, exterior domains and weighted
spaces, with $s^\ast$ encoding the geometry-dependent decay rate of the
underlying linear evolution.

\section{Conclusion and open problems}
\label{sec:conc}

We have classified the Fujita phenomenon for the reaction--diffusion
equation~\eqref{eq:main} with source intensity $h(t)=t^\gamma\ell(t)$,
where $\ell$ is slowly varying. The classification is complete for the
classical time derivative $\alpha=1$ and the full range of stable
generators $\theta\in(0,2]$. {For the Caputo derivative
$\alpha\in(0,1)$, it is complete when $N<\theta$; for $N\ge\theta$, the
blow-up range $1<p<p_c(\gamma,\alpha)$ is established in every
dimension, and a precise conjecture is formulated.} Three features stand
out. First, the slowly varying factor never shifts the critical exponent
but decides the behavior \emph{at} it. At $\alpha=1$ this happens
through the single Karamata integral $\int_1^{\infty}t^{-1}\ell(t)\,dt$;
in particular, global existence may occur at the critical exponent, in
contrast with the autonomous theory. {For $\alpha<1$ and $N<\theta$,
the boundedness of $\ell$ plays this role.} Second, the memory of the
Caputo derivative shifts the exponent of a time-dependent source,
replacing $\gamma+1$ by $\gamma+\alpha$, even though it leaves the
autonomous threshold $1+\theta/N$ untouched. {Third, for $N>\theta$
and large $\gamma$, a purely local blow-up mechanism, driven by the
singularity of the fundamental solution, overrides the scaling
exponent.}

Several directions remain open.
\begin{itemize}
\item {Prove Conjecture~\ref{conj:frac}, that is, global existence
above $p_\ast(\gamma,\alpha)=\max\{p_c(\gamma,\alpha),1+\gamma/\alpha\}$
for $N>\theta$, and determine the behavior at
$p=p_\ast(\gamma,\alpha)$. The critical-space approach of
\cite{Zhang14,CortazarQuirosWolanski24} is a natural starting point.}
\item Determine the precise blow-up rate and profile in the presence of
slowly varying coefficients, especially at the critical exponent
$p=p_c(\gamma)$ with a divergent integral $\int_1^\infty t^{-1}\ell\,dt$.
For the classical autonomous case the self-similar blow-up profile is
well understood. For slowly varying corrections, logarithmic-time
refinements should be expected.
\item Investigate the influence of a slowly varying factor on the
\emph{second} critical exponent, in the sense of Lee and
Ni~\cite{LeeNi92}, which governs the long-time behavior for initial data
with prescribed spatial decay.
\item Extend the analysis to spatially localised sources
$h(t,x)=t^\gamma\ell(t)|x|^\sigma$ of Hardy--H\'enon type, in the spirit
of~{\cite{KiraneLaskriTatar05,FinoKirane25,Majdoub23,MajdoubTorebek26}},
and to equations with an additional forcing term, as
in~{\cite{BandleLevineZhang00,BenBelgacemMajdoub26}}.
\item Treat the corresponding problem on sectorial or conical domains,
or on exterior domains, where $s^{\,\ast}$ is modified by the boundary
geometry, thereby generalising~\cite{BandleLevine89,LevineMeier90}.
\item Investigate non-monotone slowly varying functions, such as
$\ell(t)=2+\sin(\log\log t)$. The critical-case integral condition
$\int_1^{\infty}t^{-1}\ell(t)\,dt$ still determines the dichotomy, but
finer information on the blow-up time may depend on the oscillations.
\end{itemize}

\bigskip

\section*{Declarations}
\noindent\textbf{Funding.} This research received no external funding.

\noindent\textbf{Conflict of interest.} The authors declare that they have no conflicts of interest.

\noindent\textbf{Data availability.} No data were generated or analyzed in this study.
\bibliographystyle{amsalpha}

\end{document}